\documentclass{article}
\usepackage[dvips]{graphicx}
\usepackage{amssymb,amsfonts,latexsym,amsmath,amscd}

\title{On Asymptotics of Polynomial Eigenfunctions for Exactly-Solvable Differential Operators}
\author{Tanja Bergkvist\\ Department of Mathematics, University of Stockholm,\\
S-106 91 Stockholm, Sweden\\ e-mail: tanjab@math.su.se }

\begin{document}
\maketitle
\begin{abstract}
In this paper we study the class of differential operators $T=\sum_{j=1}^{k}
Q_jD^j$ with polynomial coefficients $Q_j$ in one complex variable satisfying the condition $\deg Q_j\leq j$ with equality for at least one $j$. We show that if $\deg Q_k<k$ then the root with the largest modulus of the $n$th degree eigenpolynomial $p_n$ of $T$ tends to infinity when $n\to\infty$, as opposed to the case when $\deg Q_k=k$, which we have treated previously in \cite{BR}. Moreover we present an explicit conjecture and partial results on the growth of the largest modulus of the roots of $p_n$. Based on this conjecture we deduce the algebraic equation satisfied by the Cauchy transform of the asymptotic root measure of the appropriately scaled eigenpolynomials, for which the union of all roots is conjecturally contained in a compact set. 
\end{abstract}

\section{Introduction}
In this paper we study asymptotic properties of roots in certain families of eigenpolynomials. Namely, consider a linear differential operator
\begin{displaymath}
T=\sum_{j=1}^{k}Q_{j}D^{j},
\end{displaymath}
where $D=d/dz$ and the $Q_j$ are complex polynomials in a single variable $z$ satisfying the condition $\deg Q_j\leq j$ for all $j$, and $\deg Q_k<k$ for the leading term. Such operators will be referred to as \textit{degenerate exactly-solvable operators}, see Definition 1 below. 
In this paper we study polynomial eigenfunctions of such operators, that is polynomials satisfying 
\begin{equation}\label{eigenvalueproblem}
T(p_{n})=
\lambda_{n}p_{n}
\end{equation}
 for some value of the spectral parameter $\lambda_{n}$, where $n$ is a positive integer and $\deg p_n=n$.\\

The basic motivation for this study comes from two sources: 1) a classical question going back to S. Bochner, and 2) the \textit{generalized Bochner problem}, which we describe below.\\\\
1) In 1929 Bochner asked about the classification of differential equations 
(\ref{eigenvalueproblem}) 
having an infinite sequence of \textit{orthogonal} polynomial solutions, see \cite{LJ}.
Such a system of polynomials $\{p_n\}_{n=0}^{\infty}$ which are both eigenpolynomials of some finite order differential operator \textit{and} orthogonal with respect to some suitable inner product, are referred to as \textit{Bochner-Krall orthogonal polynomial systems} (BKS), and the corresponding operators are called  \textit{Bochner-Krall operators}. It is an open problem to classify all BKS - a complete classification is only known for Bochner-Krall operators of order $k\leq 4$, and the corresponding BKS are various classical systems such as the Jacobi type, the Laguerre type, the Legendre type and the Bessel and Hermite polynomials, see \cite{EKLW}.\\\\
2) The problem of a general classisfication of linear differential operators for which the eigenvalue problem (\ref{eigenvalueproblem}) has a certain number of eigenfunctions in the form of a finite-order polynomial in some variables, is referred to as 
 \textit{the generalized Bochner problem}, see \cite{AT1} and \cite{AT2}. 
In the former paper a classification of operators possessing infinitely many finite-dimensional subspaces with a basis in polynomials is presented, and in the 
latter paper a general method has been formulated for generating eigenvalue problems for linear differential operators in one and several variables possessing polynomial solutions.\\

Notice that for the operators considered here the sequence of eigenpolynomials is in general not an orthogonal system and it can therefore not be studied by means of the extensive theory known for such systems.\\\\
\textbf{Definition 1.} We call a linear differential operator $T$ of the $k$th order \textit{exactly-solvable} if it preserves the infinite flag $\mathcal{P}_0\subset
\mathcal{P}_1\subset\mathcal{P}_2\subset\cdots\subset\mathcal{P}_n\subset\cdots $, where $\mathcal{P}_n$ is the linear space of all polynomials of degree less than or equal to $n$.\footnote{Correspondingly, a linear differential operator of the $k$th order is called \textit{quasi-exactly-solvable} if it preserves the space $\mathcal{P}_n$ for some fixed $n$.} Or, equivalently, the problem (\ref{eigenvalueproblem}) has an infinite sequence of polynomial eigenfunctions \textit{if and only if} the operator $T$ is \textit{exactly-solvable}, see \cite{AT3}.\\

Notice that any exactly-solvable operator is of the form $T=\sum_{j=1}^{k}Q_jD^j$. They split into two major classes: \textit{non-degenerate} and \textit{degenerate}, where in the former case 
$\deg Q_k=k$, and in the latter case $\deg Q_k<k$ for the leading term. 
The major difference between these two classes is that in the non-degenerate case the union of all roots of all eigenpolynomials of $T$ is contained in a compact set (see \cite{BR}), contrary to the degenerate case, which we prove in this paper.\\
 
The importance of studying eigenpolynomials of exactly-solvable operators is among other things motivated by numerous examples coming from classical orthogonal polynomials. Our study can be considered as a natural generalization of the behaviour of the maximal root for classical orthogonal polynomial families such as the Laguerre and Hermite polynomials, which appear as solutions to the eigenvalueproblem (\ref{eigenvalueproblem}) for certain choices on the polynomial coefficients $Q_j$ for a second-order degenerate exactly-solvable operator; the Laguerre polynomials appear as solutions to the differential equation $zy''(z)+(1-z)y'(z)+ny(z)=0$, and the Hermite polynomials are solutions to the differential equation $y''(z)-2zy'(z)+2ny(z)=0$ where $n$ is a nonnegative integer. Recent studies and interesting results on the asymptotic zero behaviour for these polynomials and the corresponding generalized polynomials can be found in e.g \cite{DS}, \cite{G},\cite{KMcL}, \cite{FGZ}, \cite{MMO}, \cite{KI1}, \cite{KI2} and references therein. In \cite{KI1} one can find bounds on the spacing of zeros of certain functions belonging to the Laguerre-Polya class satisfying a second order differential equation, and as a corollary new sharp inequalities on the extreme zeros of the Hermite, Laguerre and Jacobi polynomials are established.\\

Let us briefly recall our previous results.
In \cite{BR} we treated the asymptotic zero distribution for polynomial families appearing as solutions to (\ref{eigenvalueproblem}) where $T$ is an arbitrary non-degenerate exactly-solvable operator. This seems to be a natural generalization to higher orders of the Gauss hypergeometric equation. As a special case, the Jacobi polynomials appear as solutions to $(z^2-1)y''(z)+(az+b)y'(z)+cy(z)=0$, where $a,b$ and $c$ are constants satisfying 
$a>b, a+b>0$ and $c=n(1-a-n)$ for some nonnegative integer $n$. It is a classical fact that the zeros of the Jacobi polynomials lie in the interval $[-1,1]$ and that their density in this interval is proportional to $1/\sqrt{1-|z|^2}$ when the degree $n$ tends to infinity, which follows from the general theory of orthogonal polynomial systems. However, for higher-order operators of this kind, the sequence of eigenpolynomials is in general not an orthogonal system.
In \cite{BR} we proved that when $n\to\infty$, the roots of the $n$th degree eigenpolynomial $p_{n}$ for a non-degenerate exactly-solvable operator are distributed according to a certain probability measure which has compact support and which depends \textit{only} on the leading polynomial $Q_{k}$. Namely,\\\\
\textbf{Theorem A.} \textit{Let $Q_k$ be a monic polynomial of degree $k$. Then there exists a unique probability measure $\mu_{Q_k}$ with compact support whose Cauchy transform $C(z)=\int\frac{d\mu_{Q_k}(\zeta)}{z-\zeta}$ satisfies $C(z)^k=1/Q_k(z)$ for almost all $z\in\mathbb{C}$.}\\\\
\textbf{Theorem B.} \textit{Let $Q_k$ and $\mu_{Q_k}$ be as in Theorem A. Then} supp $\mu_{Q_k}$ \textit{is the union of finitely many smooth curve segments, and each of these curves is mapped to a straight line by the locally defined mapping $\Psi (z)=\int Q_k(z)^{-1/k}dz$. Moreover,} supp $\mu_{Q_k}$ \textit{contains all the zeros of $Q_k$, is contained in the convex hull of the zeros of $Q_k$ and is connected and has connected complement.}\\

If $p_n$ is a polynomial of degree $n$ we construct the probability measure $\mu_n$ by placing the point mass of size $\frac{1}{n}$ at each zero of $p_n$, and we call $\mu_n$ the \textit{root measure} of $p_n$. The following is our main result from \cite{BR}:\\\\ 
\textbf{Theorem C.} \textit{Let $p_n$ be the monic degree $n$ eigenpolynomial of a non-degenerate exactly-solvable operator, and let $\mu_n$ be the root measure of $p_n$. Then $\mu_n$ converges weakly to $\mu_{Q_k}$ when $n\to\infty$.}\\

 To illustrate, we show the zeros of some polynomial eigenfunctions for the non-degenerate exactly-solvable operator 
$T=Q_5D^5$, where $Q_5=(z-2+2i)(z+1-2i)(z+3+i)(z+2i)(z-2i-2)$. In the pictures below large dots represent the zeros of $Q_5$ and small dots represent the zeros of the eigenpolynomials 
$p_{50}, p_{75}$ and $p_{100}$ respectively.\\

\begin{tabular}{ccccc}
\includegraphics[width=2.7cm]{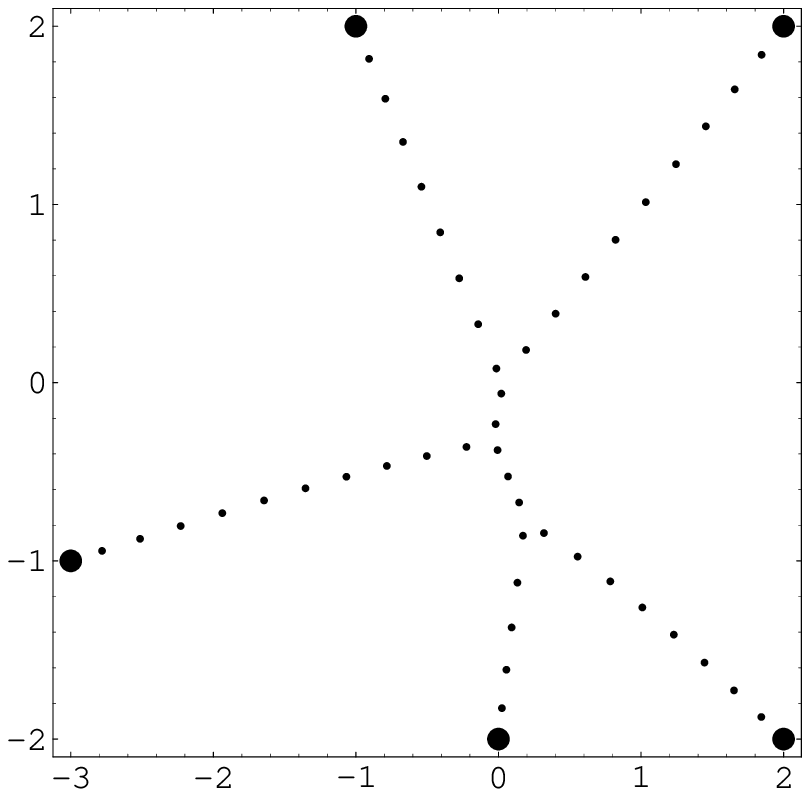} & &
\includegraphics[width=2.7cm]{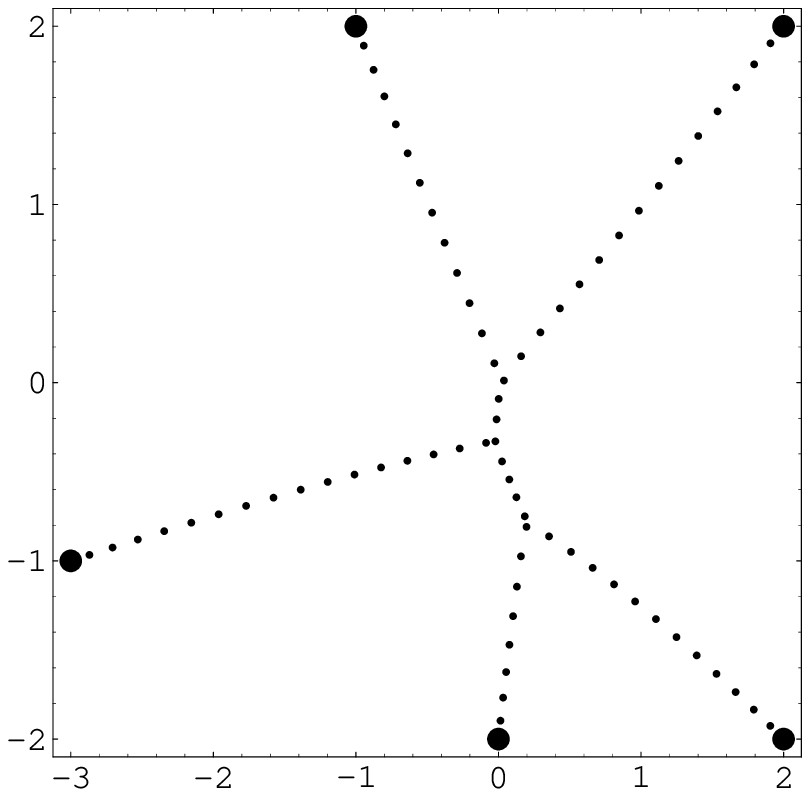} & &
\includegraphics[width=2.7cm]{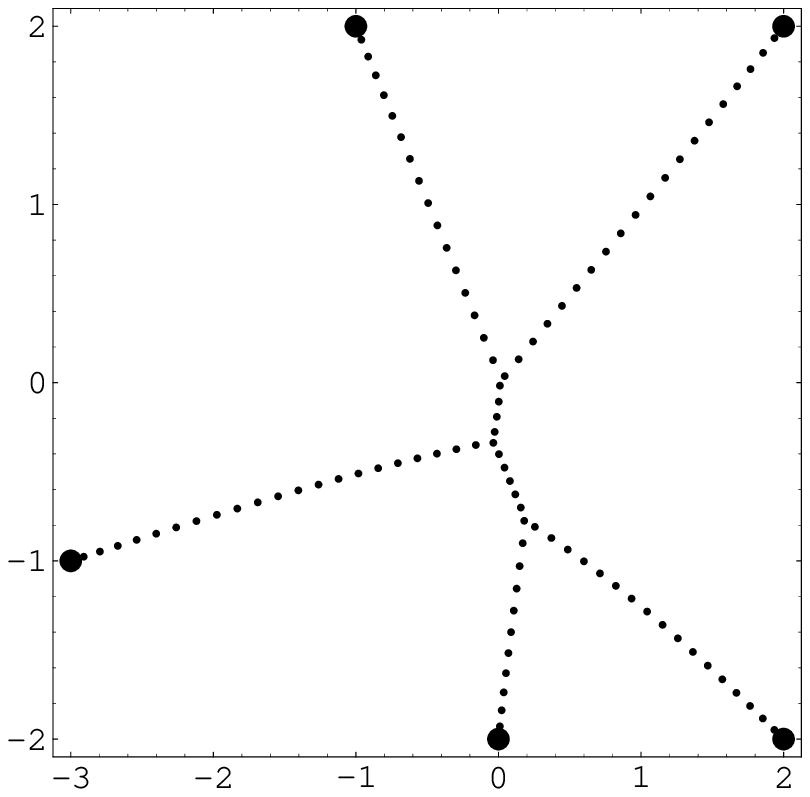}\\
n=50 & & n=75 & & n=100.
\end{tabular}\\\\

 As a consequence of the above results we were able to prove a special case of a general conjecture describing
the leading terms of all Bochner-Krall operators, see \cite{BRS}.\\

In the present paper we partially extend the above results to the case of degenerate exactly-solvable operators. Numerical evidence shows that the roots of the $n$th degree eigenpolynomial are distributed along a tree in this case too, but that the limiting root measure is compactly supported only after 
an appropraite scaling of the roots. We will assume wlog that $p_n$ is monic. 
\\\\
We start with the following preliminary result:\\\\
\textbf{Lemma 1.}
\textit{Let $T=\sum_{j=1}^{k}Q_jD^j$ be a degenerate exactly-solvable operator of order $k$. Then, for all sufficiently large integers $n$, there exists a unique constant $\lambda_n$ and a unique monic polynomial $p_n$ of degree $n$ which satisfy
 $T(p_n)=\lambda_n p_n$.}
\textit{If $\deg Q_j=j$ for {\rm precisely one} value $j<k$, 
then there exists a unique constant $\lambda_n$ and a unique monic polynomial $p_n$ of degree $n$ which
 satisfy $T(p_n)=\lambda_n p_n$ for {\rm every} integer $n=1,2,\ldots $.}\\\\
In what follows we denote by $r_n$ the largest modulus of all roots of the unique and monic $n$th degree eigenpolynomial $p_n$ of $T$, i.e. 
$r_n=\max\{|z|:p_n(z)=0\}$. These are our
\textbf{main results:}\\\\
\textbf{Theorem 1.}\footnote{This theorem is joint work with H. Rullg\aa rd.}
 \textit{Let $T$
be a degenerate exactly-solvable operator of order $k$. Then $r_n\to\infty$ when $n\to\infty$.}\\\\
Next we establish a lower bound for $r_n$ when $n\to\infty$.\\\\
 \textbf{Theorem 2.} \textit{Let $T= \sum_{j=1}^{k}Q_jD^j =\sum_{j=1}^{k}\big(\sum_{i=0}^{\deg Q_j}\alpha_{j,i}z^i\big)D^j$ be a degenerate exactly-solvable operator of order $k$. Then for any $\gamma<b$ we have}
\begin{displaymath}
\lim_{n\to\infty}\frac{r_n}{n^{\gamma}}=\infty, 
\end{displaymath}
\textit{where}
\begin{displaymath}
b:=\min_{j\in[1,k-1]}^{+}\bigg(\frac{k-j}{k-j+\deg Q_j-\deg Q_k}\bigg),
\end{displaymath}
\textit{and where the notation $\min^{+}$ means that the minimum is taken only over positive terms $(k-j+\deg Q_j-\deg Q_k)$}.\\\\
\textbf{Corollary 1.} \textit{Let $T=\sum_{j=1}^kQ_jD^j$ be a degenerate exactly-solvable operator of order $k$ such that 
$\deg Q_j\leq j_0$ for all $j>j_0$, and in particular $\deg Q_k=j_0$, where $j_0$ is the largest $j$ such that $\deg Q_j=j$.  
Then $\lim_{n\to\infty}\frac{r_n}{n^{\gamma}}=\infty$ for any $\gamma<1$.}\\\\
\textbf{Corollary 2.} \textit{Let $T=\sum_{j=1}^kQ_jD^j$ be a degenerate exactly-solvable operator of order $k$ such that 
$\deg Q_j=0$ for all $j>j_0$, where $j_0$ is the largest $j$ for which $\deg Q_j=j$. Then
$\lim_{n\to\infty}\frac{r_n}{n^{\gamma}}=\infty$ for any $\gamma<\frac{k-j_0}{k}$.}\\

In fact our extensive numerical experiments and natural heuristic arguments (see Section 3) support the following\\\\
\textbf{Main Conjecture.}
\textit{Let $T=\sum_{j=1}^{k}Q_jD^j$ be a degenerate exactly-solvable operator of order $k$ and denote by $j_0$ the largest $j$ for which $\deg Q_j=j$. Then 
\begin{displaymath}
\lim_{n\to\infty}\frac{r_n}{n^d}=c_T,
\end{displaymath}
 where $c_T>0$ is a positive constant and} 
\begin{displaymath}
d:=\max_{j\in[j_0+1,k]} \bigg(\frac{j-j_0}{j-\deg Q_j}\bigg).
\end{displaymath}\\
\textbf{Remark.} Note that Main Conjecture implies Theorem 2 since $b\leq d$.\\\\
The next two theorems support the above conjecture:\\\\
\textbf{Theorem 3.} \textit{Let $T$ be a degenerate exactly-solvable operator of order $k$ consisting of precisely two terms: $T=Q_{j_0}D^{j_0}+Q_kD^{k}$. Then there exists a positive constant $c$ such that} 
\begin{displaymath}
\lim_{n\to\infty}\inf\frac{r_n}{n^d}\geq c
\end{displaymath}
\textit{where} $d:=\max_{j\in[j_0+1,k]}
\big(\frac{j-j_0}{j-\deg Q_j}\big)=\frac{k-j_0}{k-\deg Q_k}$.\\\\
This result can be generalized to operators consisting of any number of terms, but with certain conditions on the degree of the polynomial coefficients $Q_j$ for which $j>j_0$, where $j_0$ is the largest $j$ for which $\deg Q_j=j$. Namely,\\\\
\textbf{Theorem 4.} \textit{Let $T$ be a degenerate exactly-solvable operator of order $k$. Denote 
by $j_0$ the largest $j$ such that $\deg Q_j=j$ and let
$(j-\deg Q_j)\geq (k-\deg Q_k)$ for every  $j>j_0$. Then there exists a positive constant $c>0$ such that} 
\begin{displaymath}
\lim_{n\to\infty}\inf\frac{r_n}{n^d}\geq c
\end{displaymath}
\textit{where} $d:=\max_{j\in[j_0+1,k]}
\big(\frac{j-j_0}{j-\deg Q_j}\big)=\frac{k-j_0}{k-\deg Q_k}$.\\\\

Numerical experiments show that roots of eigenpolynomials scaled according to the Main Conjecture fill certain
interesting curves in $\mathbb{C}$. To illustrate this phenomenon let us present three typical pictures. Below $p_n$ denotes the $n$th degree unique and monic eigenpolynomial of the given operator, and $q_n(z)=p_n(n^dz)$ denotes the corresponding appropriately scaled polynomial, where $d$ is as in Main Conjecture, and for which the union of all roots is (conjecturally) contained in a compact set.\\\\
\begin{tabular}{ccccc}
\textit{Fig.1:} & & \textit{Fig.2:} & & \textit{Fig.3:}\\
\includegraphics[width=2.8cm]{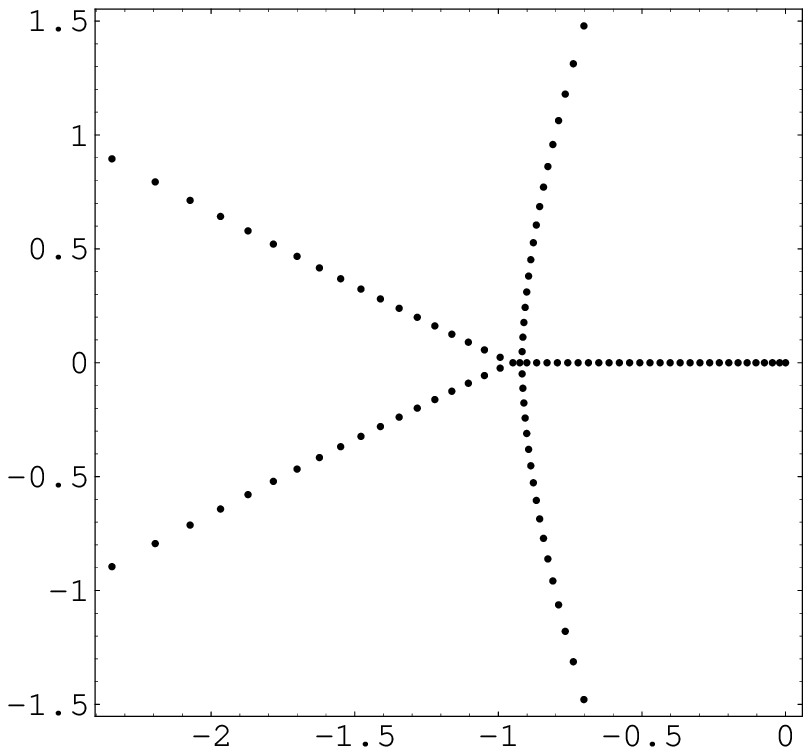} & &
\includegraphics[width=2.8cm]{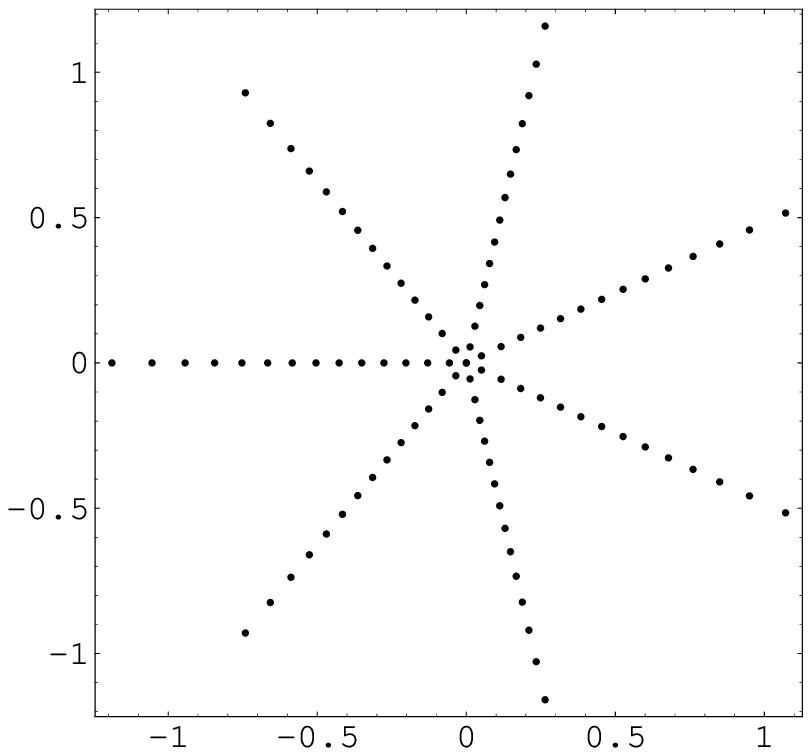} & &
\includegraphics[width=2.8cm]{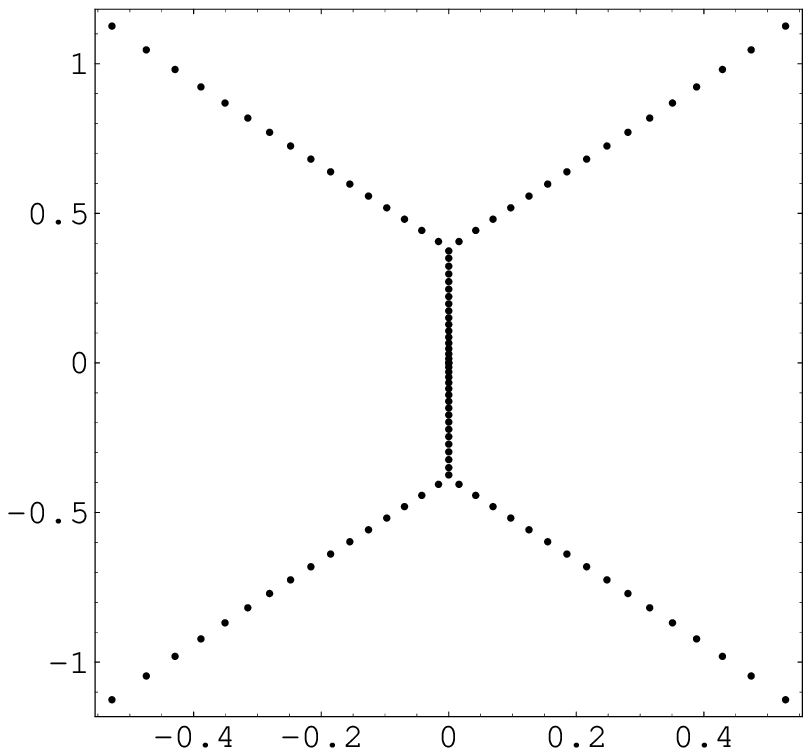}\\
roots of & & roots of & & roots of\\
 $q_{100}(z)=p_{100}(100z)$ & & $q_{100}(z)=p_{100}(100z)$
 & & $q_{100}(z)=p_{100}(100z)$ 
\end{tabular}\\\\\\
\textit{Fig.1:} $T_1=zD+zD^2+zD^3+zD^4+zD^5$.\\
\textit{Fig.2:} $T_2=z^2D^2+D^7$.\\
\textit{Fig.3:} $T_3=z^3D^3+z^2D^4+zD^5$.\\

In Section 3 we will derive the (conjectural) algebraic equation satisfied by the Cauchy transform $C(z)$ of the asymptotic root measure of the scaled eigenpolynomial
 $q_n$ for an arbitrary degenerate exactly-solvable operator. From this equation one can 
obtain detailed information about the above curves and also 
conclude which terms of the operator that are relevant for the asymptotic zero distribution of its eigenpolynomials.\footnote{As was mentioned earlier, in the non-degenerate case which we have treated previously, the asymptotic zero distribution of the eigenpolynomials depends \textit{only} on the leading coefficient $Q_k$. For the operators considered here however, the situation is more complicated.} 
\textit{Numerical evidence clearly illustrates that distinct operators whose scaled eigenpolynomials satisfy the same Cauchy transform equation when $n\to\infty$, will yield identical asymptotic zero distributions.} Below we show one such example. For further details see Section 4.3.\\\\
\begin{tabular}{ccc}
\includegraphics[width=2.8cm]{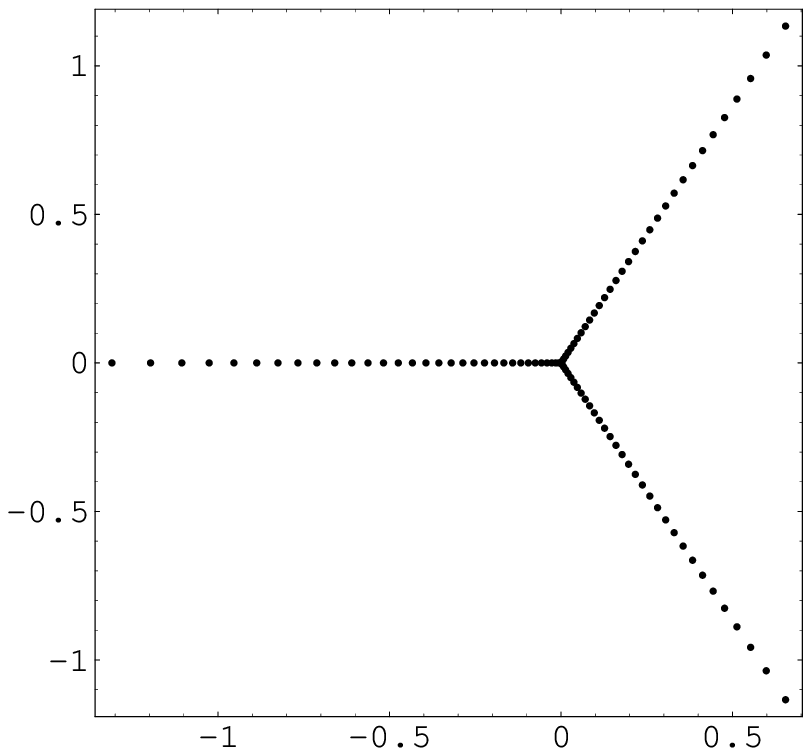} &  &
\includegraphics[width=2.8cm]{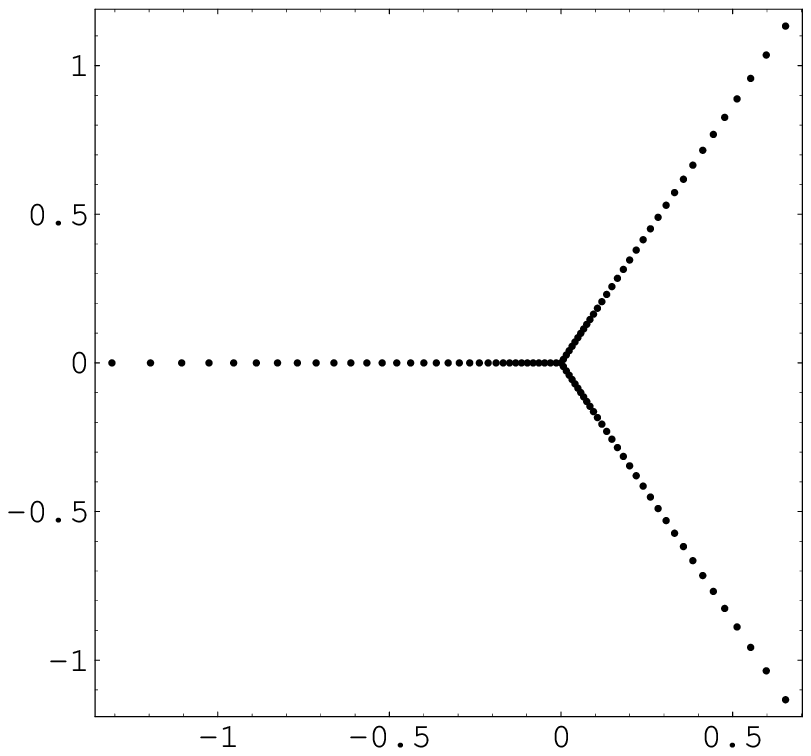}\\
$T_4=z^3D^3+z^2D^5$, & & $\widetilde{T}_4=z^2D^2+z^3D^3+zD^4+z^2D^5+D^6$,\\
roots of $q_{100}(z)=p_{100}(100^{2/3}z)$ & & roots of $q_{100}(z)=p_{100}(100^{2/3}z)$
\end{tabular}\\\\

Let us finally mention some possible applications of our results and directions for further reasearch. As was mentioned earlier, operators of the type we consider occur in the theory of Bochner-Krall orthogonal systems. A great deal is known about the asymptotic zero 
distribution of orthogonal polynomials. By comparing such known results with results on the asymptotic zero distribution of eigenpolynomials of degenerate exactly-solvable operators we believe it will be possible to gain new insight into the nature of BKS.\\

This paper is structured as follows. In Section 2 we give the proofs of the lemma, the theorems and the corollaries stated in this section. In Section 3 we explain how we arrived at Main Conjecture and how we obtained as its corollary the algebraic Cauchy transform equation. In Section 4 we give numerical evidence supporting Main Conjecture and its corollary. In Section 5 (Appendix) we give the detailed calculations led to the corollary of Main Conjecture, and we also prove that for a class of operators of the type we consider, the conjectured upper bound for $r_n$ implies the conjectured lower bound. Finally in Section 6 we discuss some open problems and directions for further research.\\\\
\textbf{Acknowledgements.} I am sincerely greatful to my PhD advisor Professor B. Shapiro for 
introducing me to this very fascinating problem and for his constant support during my work. I would also like to thank Professor J-E. Bj\"ork for stimulating discussions on the topic. I am also greatly obliged to anonymous referees for their interesting suggestions which contributed to the improvement of this paper. My research was supported by Stockholm University.    
\section{Proofs}
\textbf{Proof of Lemma 1.}
In \cite{BR} we proved that for any exactly-solvable operator $T$, the eigenvalue problem $T(p_n)=\lambda_np_n$ can be written as a linear system $MX=Y$, 
where $X$ is the coefficient vector of the monic $n$th degree eigenpolynomial $p_n$ with components 
$a_{n,0}, a_{n,1}, a_{n,2},\ldots , a_{n,n-1}$, and 
 $Y$ is a vector and $M$ is an upper triangular $n\times n$ matrix, both with entries expressible in the coefficients of $T$.
With $T=\sum_{j=1}^{k}Q_jD^j$, $Q_j=\sum_{i=0}^{\deg Q_j}\alpha_{j,i}z^i$, and 
$p_n(z)=\sum_{i=0}^{n}a_{n,i}z^i$, the eigenvalue $\lambda_n$ is given by
\begin{displaymath}
\lambda_n=\sum_{j=1}^{k}\alpha_{j,j}\frac{n!}{(n-j)!},
\end{displaymath}
and the diagonal elements of the matrix $M$ are given by
\begin{displaymath}
M_{i+1,i+1}=\sum_{1\leq j\leq\min (i,k)}^{}\alpha_{j,j}\frac{i!}{(i-j)!}-\lambda_n=\sum_{j=1}^{k}\alpha_{j,j}\bigg[\frac{i!}{(i-j)!}-\frac{n!}{(n-j)!}\bigg]
\end{displaymath}
for $i=0,1,\ldots ,n-1$. The last equality follows since 
$i!/(i-j)!=0$ for $i<j\leq k$ by definition (see Lemma 2 in \cite{BR}).
In order to prove that $p_n$ is unique we only need to check that the determinant of $M$ is nonzero, which implies that $M$ is invertible, whence the system $MX=Y$ has a unique solution. Notice that $M$ is upper triangular and thus its determinant equals the product of its diagonal elements.
We therefore prove that 
 every diagonal element $M_{i+1,i+1}$ ($i\in[0,n-1]$) is nonzero for all 
sufficiently large integers $n$ for an arbitrary $T$ as above, as well as for \textit{every} $n$ if $\deg Q_j=j$ for \textit{exactly} one $j$.\\\\ 
From the expression 
\begin{eqnarray*}
-M_{i+1,i+1}=\sum_{j=1}^{k}\alpha_{j,j}\bigg[\frac{n!}{(n-j)!}-\frac{i!}
{(i-j)!}\bigg]
\end{eqnarray*}
it is clear that $M_{i+1,i+1}\neq 0$ for every $i\in[0,n-1]$ and \textit{every} $n$ if $\alpha_{j,j}\neq 0$ for
 \textit{exactly} one $j$, that is if $\deg Q_j=j$ for precisely one $j$, and thus we have proved the second part of Lemma 1.\\

Now assume that $\deg Q_j=j$ for more than one $j$ and denote by $j_0$ the largest such $j$. Then $\alpha_{j_0,j_0}\neq 0$ and we have 
\begin{eqnarray*}
&-&M_{i+1,i+1}=\sum_{j=1}^{j_0}\alpha_{j,j}\bigg[\frac{n!}{(n-j)!}-\frac{i!}
{(i-j)!}\bigg]\\
&=&\frac{n!}{(n-j_0)!}\bigg[
\alpha_{j_0,j_0}\bigg(1-\frac{i!/(i-j_0)!}{n!/(n-j_0)!}\bigg)
+\sum_{1\leq j<j_0}^{}\alpha_{j,j}\frac{(n-j_0)!}{(n-j)!}-\sum_{1\leq j<j_0}
\frac{(n-j_0)!i!}{n!(i-j)!}\bigg].
\end{eqnarray*}
The last two sums in the brackets on the right-hand side of the above equality tend to zero when $n\to\infty$, since $j_0>j$ and $i\leq n-1$. Thus for all sufficiently large $n$ we get
\begin{eqnarray*}
-M_{i+1,i+1}=\frac{n!}{(n-j_0)!}\bigg[
\alpha_{j_0,j_0}\bigg(1-\frac{i!/(i-j_0)!}{n!/(n-j_0)!}\bigg)\bigg]\neq 0
\end{eqnarray*}
for every $i\in[0,n-1]$, and we have proved the first part of Lemma 1.
\hfill$\square$\\

To prove Theorem 1 we need the following. If $q_n$ is a polynomial of degree $n$ we construct the probability measure $\mu_n$ by placing a point mass of size $\frac{1}{n}$ at each zero of $q_n$. We call $\mu_n$ the \textit{root measure} of $q_n$.
By definition, for any polynomial $q_n$ of degree $n$, the Cauchy transform $C_{n,j}$ of the root measure $\mu^{(j)}_n$ for the $j\textit{th}$ derivative $q^{(j)}_n$ is defined by
\begin{displaymath}
C_{n,j}(z):=
\frac{q^{(j+1)}_n(z)}{(n-j)q^{(j)}_n(z)}=\int_{}^{} \frac{d\mu^{(j)}_n(\zeta)}{z-\zeta}, 
\end{displaymath}
for $j=0,1,\ldots ,n-1$,
and it is well-known that the measure $\mu_n^{(j)}$ can be reconstructed from $C_{n,j}$ by the formula $\mu_n^{(j)}=\frac{1}{\pi}\cdot\frac{\partial C_{n,j}}{\partial\bar{z}}$ where 
$\partial/\partial\bar{z}=
\frac{1}{2}(\partial/\partial x+i\partial/\partial y)$.\\\\
\textbf{Proof of Theorem 1.} Let $T=\sum_{j=1}^kQ_jD^j$ and denote by $j_0$ the largest $j$ for which 
$\deg Q_j=j$. Note that since $T$ is degenerate we have $j_0<k$. 
From the definition $C_{n,j}(z)=\frac{p_n^{(j+1)}(z)}{(n-j)p_n^{(j)}(z)}$ 
of the Cauchy transform of $p_n^{(j)}(z)$ we get 
\begin{eqnarray*}
\frac{p_n^{(j)}(z)}{p_n(z)}&=& C_{n,0}(z)C_{n,1}(z)\cdots C_{n,j-1}(z)\cdot
 n(n-1)\cdots (n-j+1)\\
&=&\frac{n!}{(n-j)!}\prod_{m=0}^{j-1}C_{n,m}(z).
\end{eqnarray*}
 With the notation $Q_j(z)=\sum_{i=0}^{\deg Q_j}\alpha_{j,i}z^i$ we have $\lambda_n=
\sum_{j=1}^{j_0}\alpha_{j,j}\frac{n!}{(n-j)!}$. Now dividing the equation 
$T(p_n)=\lambda_n p_n$ by $p_n$ we obtain
\begin{eqnarray*}
Q_k(z)\frac{p_n^{(k)}(z)}{p_n(z)}+Q_{k-1}(z)\frac{p_n^{(k-1)}(z)}{p_n(z)}+\ldots
 +Q_1(z)\frac{p_n^{\prime}(z)}{p_n(z)}
=\sum_{j=1}^{j_0}\alpha_{j,j}\frac{n!}{(n-j)!}
\end{eqnarray*}
\begin{displaymath}
\Leftrightarrow
\end{displaymath}
\begin{eqnarray}\label{egenekva}
&{}& Q_k(z)\frac{n!}{(n-k)!}\prod_{m=0}^{k-1}C_{n,m}(z)+
Q_{k-1}(z)\frac{n!}{(n-k+1)!}\prod_{m=0}^{k-2}C_{n,m}(z)+\ldots\nonumber\\
&{}&\ldots +Q_1(z)\frac{n!}{(n-1)!}C_{n,0}(z)=\sum_{j=1}^{j_0}\alpha_{j,j}\frac{n!}{(n-j)!}.
\end{eqnarray}

Dividing both sides of this equation by $\frac{n!}{(n-k)!}$ we get 
\begin{eqnarray}\label{egenekvadiv}
&{}&Q_k(z)\prod_{m=0}^{k-1}C_{n,m}(z)
\bigg[
1+\frac{(n-k)!}{(n-k+1)!}\frac{1}{C_{n,k-1}(z)}\frac{Q_{k-1}(z)}{Q_k(z)}+\nonumber\\
&{}&\frac{(n-k)!}{(n-k+2)!}\frac{1}
{C_{n,k-1}(z)C_{n,k-2}(z)}\frac{Q_{k-2}(z)}{Q_k(z)}+ \ldots\nonumber\\ &{}&\ldots +\frac{(n-k)!}{(n-1)!}\frac{1}
{\prod_{m=1}^{k-1}C_{n,m}(z)}
\frac{Q_1(z)}{Q_k(z)}\bigg]
=\sum_{j=1}^{j_0}\alpha_{j,j}\frac{(n-k)!}{(n-j)!}.
\end{eqnarray}
Now assume that all zeros of all $p_n$ are uniformly bounded. Then we can take a subsequence $\{p_{n_i}\}$ such that all the corresponding root measures 
$\mu_{n_i}$ are weakly convergent to a compactly supported probability measure. Then all Cauchy transforms $C_{n_i,m}$ will be uniformly convergent to a non-vanishing function outside some large disc, which in particular contains all the roots of $Q_k(z)$. Since $j_0<k$, the right-hand side of (\ref{egenekvadiv}) tends to zero when $n\to\infty$.  On the other hand, in the left-hand side of (\ref{egenekvadiv}), all  
terms in the bracket except for the constant term $1$ tend to zero when $n\to\infty$, and thus the limit of the left-hand side equals
$\lim_{n\to\infty}Q_k(z)\prod_{m=0}^{k-1}C_{n,m}(z)=K\neq 0$, and we obtain a contradiction when $n\to\infty$. 
\hfill$\square$\\\\
In order to prove Theorem 2 we need the following two lemmas, where Lemma 2 is used to prove Lemma 3. \\\\
\textbf{Lemma 2.} \textit{Let $z_n$ be a root of $p_n$ with the largest modulus $r_n$. Then, for any complex number $z_0$ such that $|z_0|=r_0\geq r_n$, we have $|C_{n,j}(z_0)|\geq\frac{1}{2r_0}$ for all $j\geq 0$.}\\\\
\textbf{Proof.} Recall that $C_{n,j}(z):=\int_{}^{} \frac{d\mu^{(j)}_n(\zeta)}{z-\zeta}=\frac{p^{(j+1)}_n(z)}{(n-j)p^{(j)}_n(z)}$.
With $\zeta$ being some root of $p^{(j)}_n(z)$ we have $|\zeta |\leq |z_0|$ by Gauss-Lucas theorem. Thus 
$\frac{1}{z_0-\zeta}=\frac{1}{z_0}\cdot\frac{1}{1-\zeta /z_0}=\frac{1}{z_0}\cdot\frac{1}{1-\theta}$ where $|\theta|=
|\zeta/z_0|\leq 1$. With $w=\frac{1}{1-\theta}$ we obtain
\begin{displaymath}
|w-1|=\frac{|\theta|}{|1-\theta|}=|\theta||w|\leq |w|\Leftrightarrow |w-1|\leq |w|\Rightarrow Re(w)\geq 1/2,
\end{displaymath}
and thus
\begin{eqnarray*}
|C_{n,j}(z_0)| & = & \bigg|\int\frac{d\mu^{(j)}_n(\zeta)}{z_0-\zeta}\bigg|=\frac{1}{r_0}\bigg|\int\frac{d\mu^{(j)}_n(\zeta)}{1-\theta}\bigg|=\frac{1}{r_0}\bigg|\int wd\mu^{(j)}_n(\zeta)\bigg|\\
& \geq & \frac{1}{r_0}\bigg|\int Re(w)d\mu^{(j)}_n(\zeta)\bigg|\geq\frac{1}{2r_0}\int d\mu^{(j)}_n(\zeta)=\frac{1}{2r_0}.
\end{eqnarray*}
\hfill$\square$\\\\
\textbf{Lemma 3.}
\textit{Let 
$T=\sum_{j=1}^{k}Q_jD^j= \sum_{j=1}^{k}\big(\sum_{i=0}^{\deg Q_j}\alpha_{j,i}z^i\big)D^j$
be a degenerate exactly-solvable operator of order $k$. Wlog we assume that $Q_k$ is monic, i.e. $\alpha_{k,\deg Q_k}=1$.
 Let $z_n$ be a root of $p_n$ with the largest modulus $r_n$. Then the following inequality holds:}
\begin{equation}\label{lemma3}
1\leq \sum_{j=1}^{k-1}\sum_{i=0}^{\deg Q_j}|\alpha_{j,i}|2^{k-j}\frac{r_n^{k-j-\deg Q_k+i}}{(n-k+1)^{k-j}}+\sum_{0\leq i<\deg Q_k}\frac{|\alpha_{k,i}|}{r_{n}^{\deg Q_k-i}}.
\end{equation}
\textbf{Proof.}
From $C_{n,j}(z)=\frac{p^{(j+1)}_n(z)}{(n-j)p^{(j)}(z)}$ we get
\begin{eqnarray}\label{rel}
p^{(j)}(z)=
\frac{p^{(k)}_n(z)}{(n-k+1)(n-k+2)\cdots (n-j)\prod_{m=j}^{k-1}C_{n,m}(z)}\quad \forall\quad j<k.
\end{eqnarray}
 Inserting $z_n$ in the eigenvalue equation $Tp_n(z)=\lambda_np_n(z)$ we obtain
\begin{displaymath}
  \sum_{j=1}^{k-1}\bigg(\sum_{i=0}^{\deg Q_j}\alpha_{j,i}z_n^i\bigg)p^{(j)}_n(z_n)+\bigg(
\sum_{i=0}^{\deg Q_k}\alpha_{k,i}z_n^i\bigg)p^{(k)}_n(z_n)=\lambda_np_n(z_n)=0,
\end{displaymath}
and after division by $z_n^{\deg Q_k}p^{(k)}_n(z_n)$ we obtain
\begin{displaymath}
\sum_{j=1}^{k-1}\bigg(\sum_{i=0}^{\deg Q_j}\alpha_{j,i}\frac{1}{z_n^{\deg Q_k-i}}\bigg)
\frac{p^{(j)}_n(z_n)}{p^{(k)}_n(z_n)}+
\sum_{0\leq i<\deg Q_k}^{}\alpha_{k,i}\frac{1}{z_n^{\deg Q_k-i}}+1=0.
\end{displaymath}
Thus, applying (\ref{rel}) and Lemma 2, we obtain  
\begin{eqnarray*}
1&=&\bigg|\sum_{j=1}^{k-1}\bigg(\sum_{i=0}^{\deg Q_j}\alpha_{j,i}\frac{1}{z_n^{\deg Q_k-i}}\bigg)
\frac{p^{(j)}_n(z_n)}{p^{(k)}_n(z_n)}+
\sum_{0\leq i<\deg Q_k}^{}\alpha_{k,i}\frac{1}{z_n^{\deg Q_k-i}}\bigg|\\
&\leq &\sum_{j=1}^{k-1}\bigg|\sum_{i=0}^{\deg Q_j}\alpha_{j,i}\frac{1}{z_n^{\deg Q_k-i}}\bigg|\frac{|p^{(j)}_n(z_n)|}{|p^{(k)}_n(z_n)|}+\sum_{0\leq i<\deg Q_k}\frac{|\alpha_{k,i}|}{r_{n}^{\deg Q_k-i}}\\
&\leq & \sum_{j=1}^{k-1}\sum_{i=0}^{\deg Q_j}\frac{|\alpha_{j,i}|}{r_n^{\deg Q_k-i}}\frac{1}{(n-k+1)\cdots (n-j)\prod_{m=j}^{k-1}|C_{n,m}(z_n)|}+\sum_{0\leq i<\deg Q_k}\frac{|\alpha_{k,i}|}{r_{n}^{\deg Q_k-i}}\\
&\leq &\sum_{j=1}^{k-1}\sum_{i=0}^{\deg Q_j}\frac{|\alpha_{j,i}|}{r_n^{\deg Q_k-i}}\frac{(2r_n)^{k-j}}{(n-k+1)^{k-j}}+\sum_{0\leq i<\deg Q_k}\frac{|\alpha_{k,i}|}{r_{n}^{\deg Q_k-i}}\\
&=&\sum_{j=1}^{k-1}\sum_{i=0}^{\deg Q_j}|\alpha_{j,i}|2^{k-j}\frac{r_n^{k-j-\deg Q_k+i}}{(n-k+1)^{k-j}}+\sum_{0\leq i<\deg Q_k}\frac{|\alpha_{k,i}|}{r_{n}^{\deg Q_k-i}}.
\end{eqnarray*}
\hfill$\square$\\\\
The proof of Theorem 2 follows from Theorem 1 and Lemma 3.\\\\
\textbf{Proof of Theorem 2.}
Applying Theorem 1, we see that the last sum on the right-hand side of 
inequality (\ref{lemma3}) in Lemma 3
tends to zero when $n\to\infty$. 
Now consider the double sum on the right-hand side of (\ref{lemma3}). If, for given $i$ and $j$, the exponent $(k-j-\deg Q_k+i)$ of $r_n$ is negative or zero, the corresponding term tends to zero when $n\to\infty$ by Theorem 1. We now consider the remaining terms in the double sum, namely those for which the exponent $(k-j-\deg Q_k+i)$ of $r_n$ is positive.
If $r_n\leq c_0(n-k+1)^{\gamma}$ where $c_0>0$ and $\gamma<\frac{k-j}{k-j+i-\deg Q_k}$ for given $j\in[1,k-1]$ and given $i\in[0,\deg Q_j]$, then the corresponding term 
\begin{displaymath}
\frac{r^{k-j+i-\deg Q_k}_n}{(n-k+1)^{k-j}}=\bigg(
\frac{r_n}{(n-k+1)^{\frac{k-j}{k-j+i-\deg Q_k}}}\bigg)^{k-j+i-\deg Q_k}
\end{displaymath}
in the double sum tends to zero when $n\to\infty$. Thus 
assume that $r_n\leq c_0(n-k+1)^{\gamma}$ where $c_0$ is a positive constant and $\gamma <b$ where 
\begin{displaymath}
b=\min_{{j\in[1,k-1]\atop i\in[0,j]}}^{+}\frac{k-j}{k-j+i-\deg Q_k}
=\min_{j\in[1,k-1]}^{+}\frac{k-j}{k-j+\deg Q_j-\deg Q_k},
\end{displaymath}
and where the notation $\min^+$ means that we only take the minimum over positive terms $(k-j+i-\deg Q_k)$ and $(k-j+\deg Q_j-\deg Q_k)$.\footnote{On the left-hand side in the expression for $b$ above we take the minimum over $i\in[0,\deg Q_j]$, so we can put $i=\deg Q_j$ in this expression. Thus with 
$b=\min_{j\in[1,k-1]}^{+}\frac{k-j}{k-j+\deg Q_j-\deg Q_k}$ 
 we get that 
$\gamma<\frac{k-j}{k-j+i-\deg Q_k}$ for \textit{every} $j\in[1,k-1]$ and \textit{every} $i\in[0,\deg Q_j]$. Then if $r_n\leq c_0(n-k+1)^{\gamma}$ and $\gamma<b$, every term with positive exponent $(k-j+i-\deg Q_k)$ will tend to zero when $n\to\infty$.} 
Then every term in the double sum tends to zero when $n\to\infty$, and we obtain a contradiction to (\ref{lemma3}) when $n\to\infty$.
Thus for all sufficiently large integers $n$ we must have $r_n>c_0(n-k+1)^{\gamma}$ for all $\gamma <b$,   
and hence $\lim\inf_{n\to\infty}\frac{r_n}{n^{\gamma}}>c_0$ for any $\gamma<b$.
But for any such $\gamma$ we can form $\gamma^{'}=\frac{\gamma +b}{2}$ for which $\gamma^{'}<b$ and $\gamma < \gamma^{'}$, and thus  
$\lim_{n\to\infty}\frac{r_n}{n^{\gamma}}=\infty$ for all $\gamma<b$.
\hfill $\square$\\\\
\textbf{Proof of Corollary 1.} For this class of operators we have 
\begin{eqnarray*}
b&:=& \min_{j\in[1,k-1]}^{+}\frac{k-j}{k-j+\deg Q_j-\deg Q_k}\\
&=& \min_{j\in[1,k-1]}^{+}\frac{k-j}{k-j+\deg Q_j-j_0}
=\frac{k-j_0}{k-j_0}=1,
\end{eqnarray*}
and the proof is complete applying Theorem 2.\hfill$\square$\\\\
\textbf{Proof of Corollary 2.}
For this class of operators we have   
\begin{eqnarray*}
b&:=&\min^+_{j\in[1,k-1]}\bigg(\frac{k-j}{k-j+\deg Q_j-\deg Q_k}\bigg)\\
&=&\min^+_{j\in[1,k-1]}\bigg(\frac{k-j}{k-j+\deg Q_j}\bigg)=
\min_{j\in[1,j_0]}\frac{k-j}{k}=\frac{k-j_0}{k}
\end{eqnarray*}
where the third equality follows from  
choosing any $j$ for which $\deg Q_j=j$, and the minimum is then attained for $j=j_0$ (note that for $j>j_0$ we get $(k-j)/(k-j+\deg Q_j)=1>(k-j_0)/k)$, and the proof is complete applying Theorem 2.
\hfill$\square$\\\\ 
\textbf{Remark}. Note that for the class of operators considered in Corollary 1
the Main Conjecture claims that $\lim_{n\to\infty}
\frac{r_n}{n}=c_{T}$ for some $c_T>0$, since  
$d:=\max_{j\in[j_0+1,k]}\big(\frac{j-j_0}{j-\deg Q_j}\big)=
\frac{k-j_0}{k-j_0}=1$ (the maximum is attained by choosing any $j>j_0$ such that $\deg Q_j=j_0$, e.g. $j=k$), and for the class of operators considered in Corollary 2 the
Main Conjecture claims that $\lim_{n\to\infty}
\frac{r_n}{n^{(k-j_0)/k}}=c_T$ for some $c_T>0$, since 
\begin{displaymath}
d:=\max_{j\in[j_0+1,k]}\bigg(\frac{j-j_0}{j-\deg Q_j}\bigg)=
\max_{j\in[j_0+1,k]}\bigg(\frac{j-j_0}{j}\bigg)=\frac{k-j_0}{k}.
\end{displaymath}\\
\textbf{Remark.} For a class of operators containing the operators considered in Corollaries 1 and 2 we can actually prove that the conjectured upper bound $\lim_{n\to\infty}\sup\frac{r_n}{n^d}\leq c_1$ implies the conjectured lower bound 
$\lim_{n\to\infty}\inf\frac{r_n}{n^d}\geq c_0$ where $c_1\geq c_0>0$, see Theorem 5 in Section 5.2.\\\\
\textbf{Proof of Theorem 3.} Clearly $\deg Q_{j_0}=j_0$ since there exists at least one such $j<k$. Set  
\begin{displaymath}
T=Q_{j_0}D^{j_0}+Q_kD^{k}=\sum_{i=0}^{j_0}
\alpha_{j_0,i}z^iD^{j_0}+\sum_{i=0}^{\deg Q_k}\alpha_{k,i}z^iD^{k},
\end{displaymath}
where $\alpha_{j_0,j_0}\neq 0$, and where we wlog assume that $Q_k$ is monic. From inequality (\ref{lemma3}) in Lemma 3 we have 
\begin{eqnarray}
1&\leq& \sum_{i=0}^{j_0}|\alpha_{j_0,i}|
2^{k-j_0}\frac{r^{i-\deg Q_k+k-j_0}_n}{(n-k+1)^{k-j_0}}
+\sum_{0\leq i<\deg Q_k}^{}|\alpha_{k,i}|\frac{1}{r^{\deg Q_k-i}_n}\nonumber\\
&\leq & \sum_{i=0}^{j_0}|\alpha_{j_0,i}|
2^{k-j_0}\frac{r^{i-\deg Q_k+k-j_0}_n}{(n-k+1)^{k-j_0}}
+\epsilon,\nonumber
\end{eqnarray}
where we choose $n$ so large that $\epsilon<1$ (this is possible since $\epsilon\to 0$ when $n\to\infty$ due to Theorem 1). Thus for sufficiently large $n$ we get
\begin{eqnarray*}
c_0&\leq& \sum_{i=0}^{j_0}|\alpha_{j_0,i}|
2^{k-j_0}\frac{r^{i-\deg Q_k+k-j_0}_n}{(n-k+1)^{k-j_0}}\\
&\leq& 
\sum_{i=0}^{j_0}|\alpha_{j_0,i}|2^{k-j_0}
\frac{r^{k-\deg Q_k}_n}{(n-k+1)^{k-j_0}}\\
&=& K\frac{r^{k-\deg Q_k}_n}{(n-k+1)^{k-j_0}},
\end{eqnarray*}
where $1-\epsilon=c_0\to 1$ when $n\to\infty$, and $K>0$ since $\alpha_{j_0,j_0}\neq 0$ (the second inequality follows since $i\leq j_0$). 
 Thus 
\begin{displaymath}
r_n\geq \big(\frac{c_0}{K}\big)^{1/(k-\deg Q_k)}(n-k+1)^{\frac{k-j_0}{k-\deg Q_k}}
\end{displaymath}
for sufficiently large integers $n$, and hence there exists a positive constant $c=(1/K)^{1/(k-\deg Q_k)}$ such that 
\begin{displaymath}
\lim_{n\to\infty}\inf\frac{r_n}
{n^{\big(\frac{k-j_0}{k-\deg Q_k}\big)}}\geq c.
\end{displaymath}
Finally, it is clear that for this two-term operator 
\begin{displaymath}
d:=\max_{j\in[j_0+1,k]}
\bigg(\frac{j-j_0}{j-\deg Q_j}\bigg)=\frac{k-j_0}{k-\deg Q_k},
\end{displaymath}
and we are done.
\hfill$\square$\\\\
\textbf{Remark.} If, in Theorem 3, $Q_k$ is a monomial (i.e. $Q_k=z^{\deg Q_k}$), then there exists a positive constant $c$ such that $r_n\geq c(n-k+1)^{d}$ for \textit{every} $n$, where $d:=\max_{j\in[j_0+1,k]}\big(\frac{j-j_0}
{j-\deg Q_j}\big)=\frac{k-j_0}{k-\deg Q_k}$.
 This is easily seen from the calculations in the proof of Theorem 3 since $\sum_{0\leq i<\deg Q_k}^{}|\alpha_{k,i}|\frac{1}{r^{\deg Q_k-i}_n}$ on the right-hand side of (\ref{lemma3}) vanishes, and therefore $1\leq
 K\frac{r^{k-\deg Q_k}_n}{(n-k+1)^{k-j_0}}$ for \textit{every} $n$. From the second part of Lemma 1 we know that for this class of operators there exists a unique eigenpolynomial $p_n$ for every $n$, and the conclusion follows.\\\\
\textbf{Proof of Theorem 4.} For this class of operators $(j-\deg Q_j)\geq(k-\deg Q_k)$ for every $j>j_0$ and thus 
\begin{displaymath}
d:=\max_{j\in[j_0+1,k]}\bigg(\frac{j-j_0}{j-\deg Q_j}\bigg)=\frac{k-j_0}{k-\deg Q_k}.
\end{displaymath} 
Assuming that $Q_k$ is monic we have the 
inequality
\begin{eqnarray}\label{ineqThm4}
1\leq \sum_{j=1}^{k-1}\sum_{i=0}^{\deg Q_j}|\alpha_{j,i}|2^{k-j}\frac{r_n^{k-j+i-\deg Q_k}}
{(n-k+1)^{k-j}}+\sum_{0\leq i<\deg Q_k}\frac{|\alpha_{k,i}|}{r_{n}^{\deg Q_k-i}}
\end{eqnarray}
by Lemma 3. 
The last sum here tends to zero when $n\to\infty$ by Theorem 1.
Considering the double sum on the right-hand side of (\ref{ineqThm4}) we see that for every $j$ we have (since
$i\leq \deg Q_j$) that  
\begin{eqnarray}\label{equality}
&{}&\sum_{i=0}^{\deg Q_j}|\alpha_{j,i}|2^{k-j}\frac{r_n^{k-j+i-\deg Q_k}}
{(n-k+1)^{k-j}}=
\sum_{i=0}^{\deg Q_j}|\alpha_{j,i}|2^{k-j}\frac{r_n^{k-j+\deg Q_j-\deg Q_k}}
{(n-k+1)^{k-j}}r_n^{i-\deg Q_j}\nonumber\\
&=& \frac{r_n^{k-j+\deg Q_j-\deg Q_k}}{(n-k+1)^{k-j}}
\bigg(2^{k-j}|\alpha_{j,\deg Q_j}|+\sum_{i<\deg Q_j}2^{k-j}|\alpha_{j,i}|r_n^{i-\deg Q_j}\bigg)\nonumber\\
&=& K_{j,n}\frac{r_n^{k-j+\deg Q_j-\deg Q_k}}{(n-k+1)^{k-j}},
\end{eqnarray}
where
\begin{displaymath}
K_{j,n}=2^{k-j}|\alpha_{j,\deg Q_j}|+\sum_{i<\deg Q_j}2^{k-j}|\alpha_{j,i}|r_n^{i-\deg Q_j},
\end{displaymath}
where $K_{j,n}>0$ 
 since $\alpha_{j,\deg Q_j}\neq 0$. Also, $K_{j,n}<\infty$ when $n\to\infty$ since 
$(i-\deg Q_j)<0$ for the exponent of $r_n$, and then using Theorem 1 (note that $K_{j,n}\to 2^{k-j}|\alpha_{j,\deg Q_j}|$ 
when $n\to\infty$).
With the decomposition\\
 $A=\{j:\deg Q_j=j\}$,\\
 $B=\{j:\deg Q_j<j\quad \textit{and}\quad (k-j+\deg Q_j-\deg Q_k)>0\}$,\\
 $C=\{j:\deg Q_j<j\quad \textit{and}\quad (k-j+\deg Q_j-\deg Q_k)\leq 0\},$\\
 and using (\ref{equality}) we see that inequality (\ref{ineqThm4}) is equivalent to:
\begin{eqnarray*}
1&\leq& \sum_{j=1}^{k-1}\sum_{i=0}^{j}|\alpha_{j,i}|2^{k-j}\frac{r_n^{k-j+i-\deg Q_k}}{(n-k+1)^{k-j}}+\sum_{0\leq i<\deg Q_k}\frac{|\alpha_{k,i}|}{r_{n}^{\deg Q_k-i}}\\
&= & \sum_{j\in A}^{}K_{j,n}\frac{r_n^{k-\deg Q_k}}{(n-k+1)^{k-j}}+
\sum_{j\in B}^{}K_{j,n}\frac{r_n^{k-j+\deg Q_j-\deg Q_k}}{(n-k+1)^{k-j}}\\
&+& \sum_{j\in C}^{}K_{j,n}\frac{r_n^{k-j+\deg Q_j-\deg Q_k}}{(n-k+1)^{k-j}}
+\sum_{0\leq i<\deg Q_k}^{}\frac{|\alpha_{k,i}|}{r_{n}^{\deg Q_k-i}}.\\
\end{eqnarray*}

The last two sums on the right-hand side of this inequality both tend to zero when $n\to\infty$, the last one due to Theorem 1, and the sum over $C$ since $(j-\deg Q_j)\geq (k-\deg Q_k)$ $\Leftrightarrow$ $(k-j+\deg Q_j-\deg Q_k)\leq 0$ for every $j\in C$ by assumption, and then applying Theorem 1.
Therefore, when $n\to\infty$, we get the inequality 
\begin{eqnarray}\label{basicineq}
c_0\leq \sum_{j\in A}^{}K_{j,n}\frac{r_n^{k-\deg Q_k}}{(n-k+1)^{k-j}}+
\sum_{j\in B}^{}K_{j,n}\frac{r_n^{k-j+\deg Q_j-\deg Q_k}}{(n-k+1)^{k-j}}
\end{eqnarray}
where
\begin{displaymath}
c_0=1-\sum_{j\in C}^{}K_{j,n}\frac{r_n^{k-j+\deg Q_j-\deg Q_k}}{(n-k+1)^{k-j}}-
\sum_{0\leq i<\deg Q_k}\frac{|\alpha_{k,i}|}{r_{n}^{\deg Q_k-i}}.
\end{displaymath}
Note that $c_0\to 1$ when $n\to\infty$.

Now assume that $B$ is empty. This corresponds to an operator such that $(j-\deg Q_j)\geq (k-\deg Q_k)$ for every $j$ for which $\deg Q_j<j$. Then inequality (\ref{basicineq}) above becomes 
\begin{eqnarray}\label{sumoverA}
c_0&\leq &\sum_{j\in A}^{}K_{j,n}\frac{r_n^{k-\deg Q_k}}{(n-k+1)^{k-j}}\nonumber \\
&=&
\frac{r_n^{k-\deg Q_k}}{(n-k+1)^{k-j_0}}\bigg(
K_{j_0,n}+\sum_{j\in A\backslash \{j_0\}}K_{j,n}\frac{1}{(n-k+1)^{j_0-j}}\bigg)\nonumber \\
&\leq & K_A\frac{r_n^{k-\deg Q_k}}{(n-k+1)^{k-j_0}}
\end{eqnarray}
where $K_A$ is a positive constant which is finite when $n\to\infty$, since 
$j_0-j>0$ for every $j\in A\backslash \{j_0\}$ (recall that $j_0$ is the largest element in $A$ by definition).  
Thus for all sufficiently large integers $n$ we have 
\begin{displaymath}
r_n\geq \big(\frac{c_0}{K_A}\big)^{1/(k-\deg Q_k)}(n-k+1)^
{\frac{k-j_0}{k-\deg Q_k}}, 
\end{displaymath}
and therefore there exists a positive constant $c=(1/K_A)^{1/(k-\deg Q_k)}$ such that
\begin{eqnarray*}
\lim_{n\to\infty}\inf\frac{r_n}{n^{\big(\frac{k-j_0}{k-\deg Q_k}\big)}}\geq c,
\end{eqnarray*}
and we are done.

Now assume that $B$ is nonempty. Again inequality (\ref{basicineq}) holds, i.e. 
\begin{eqnarray*}
c_0&\leq & \sum_{j\in A}^{}K_{j,n}\frac{r_n^{k-\deg Q_k}}{(n-k+1)^{k-j}}+
\sum_{j\in B}^{}K_{j,n}\frac{r_n^{k-j+\deg Q_j-\deg Q_k}}{(n-k+1)^{k-j}},
\end{eqnarray*}
where $c_0\to 1$ when $n\to\infty$. 
From (\ref{sumoverA}) we have 
\begin{displaymath}
\sum_{j\in A}^{}K_{j,n}\frac{r_n^{k-\deg Q_k}}{(n-k+1)^{k-j}}
\leq K_A\frac{r_n^{k-\deg Q_k}}{(n-k+1)^{k-j_0}}
\end{displaymath}
for the sum over $A$ 
for large $n$ and thus 
\begin{eqnarray*}
c_0&\leq & \sum_{j\in A}^{}K_{j,n}\frac{r_n^{k-\deg Q_k}}{(n-k+1)^{k-j}}+
\sum_{j\in B}^{}K_{j,n}\frac{r_n^{k-j+\deg Q_j-\deg Q_k}}{(n-k+1)^{k-j}}\\
&\leq & K_A\frac{r_n^{k-\deg Q_k}}{(n-k+1)^{k-j_0}}+\sum_{j\in B}^{}K_{j,n}\frac{r_n^{k-j+\deg Q_j-\deg Q_k}}{(n-k+1)^{k-j}}\\
&= &\frac{r_n^{k-\deg Q_k}}{(n-k+1)^{k-j_0}}\bigg(K_A+
\sum_{j\in B}K_{j,n}\frac{r_n^{\deg Q_j-j}}{(n-k+1)^{j_0-j}}\bigg)\\
&\leq & K_{AB}\frac{r_n^{k-\deg Q_k}}{(n-k+1)^{k-j_0}},
\end{eqnarray*}
where $K_{AB}$ is a positive and finite constant when $n\to\infty$ (note that $K_{AB}\to K_A$ when $n\to\infty$, since
$(\deg Q_j-j)<0$ and $j_0-j>0$ for every $j\in B$).
Thus for all sufficiently large integers $n$ we have 
\begin{displaymath}
r_n\geq \big(\frac{c_0}{K_{AB}}\big)^{1/(k-\deg Q_k)}(n-k+1)^
{\frac{k-j_0}{k-\deg Q_k}},
\end{displaymath}
so there exists a positive constant $c=(1/K_{AB})^{1/(k-\deg Q_k)}$ such that 
\begin{eqnarray*}
\lim_{n\to\infty}\inf\frac{r_n}{n^{\big(\frac{k-j_0}{k-\deg Q_k}\big)}}\geq c.
\end{eqnarray*}
\hfill$\square$
\section{Main Conjecture and its Corollary} 
In this section we explain how we arrived at Main Conjecture (see Section 1)
and obtain as a corollary of our method the (conjectural) algebraic equation satisfied by the Cauchy transform of the asymptotic root measure of the properly scaled eigenpolynomials.\\

\qquad\qquad\qquad\textbf{How did we arrive at Main Conjecture?}\\

Let $T=\sum_{j=1}^{k}Q_jD^j=\sum_{j=1}^{k}\big(\sum_{i=0}^{\deg Q_j}\alpha_{j,i}z^i\big)D^j$ be an arbitrary degenerate exactly-solvable operator of order $k$ and denote by $j_0$ the largest $j$ for which $\deg Q_j=j$. Wlog we assume that $Q_{j_0}$ is monic, i.e. $\alpha_{j_0,j_0}=1$. Consider the \textit{scaled eigenpolynomial}
$q_n(z)=p_n(n^{d}z)$, where $p_n(z)$ is the unique and monic $n$th degree eigenpolynomial of $T$, and $d$ is some real number. The goal is now to obtain a well-defined algebraic equation for the Cauchy transform of the root measure $\mu_n$ of the scaled eigenpolynomial $q_n$ when $n\to\infty$, and as we will see in the process of doing this, we are forced to choose $d$ as in Main Conjecture.\footnote{It is already well-known that for the Laguerre polynomials, which appear as eigenpolynomials for a second order exactly-solvable operator, the largest root grows as $n$ when $n\to\infty$ and thus $d=1$ in this case, which is consistent with Main Conjecture.}\\\\
\textbf{Basic assumption.} When performing our calculations we assume that the root measures  $\mu^{(0)}_n, \mu^{(1)}_n, \mu^{(2)}_n\ldots, \mu^{(k-1)}_n$ of the \textit{scaled eigenpolynomial} $q_n(z)$ and its derivatives up to the $k$th order exist when $n\to\infty$ and that they are are all weakly convergent to the \textit{same} asymptotic root measure $\mu$.\footnote{Conjecturally supp $\mu$ is a tree, see Section 6 on Open Problems.} Thus the corresponding Cauchy transforms are all asymptotically identical, and we define
$C(z):=\lim_{n\to\infty}C_{n,j}(z)$ for all $j\in[0,k-1]$, where $C(z)$ is the Cauchy transform of $\mu$ and is considered for $z$'s away from the support of $\mu$. Computer experiments strongly indicate that this assumption is true - for details see Section 4.2.\\

From the definition of the Cauchy transform we obtain
\begin{eqnarray}
  \prod_{i=0}^{j-1}C_{n,i}(z)&=&\prod_{i=0}^{j-1}\frac{q^{(i+1)}_n(z)}{(n-j)q^{(i)}_n(z)}\nonumber\\
&=&\frac{q^{(1)}_n(z)}{nq_{n}(z)}\cdot\frac{q^{(2)}_n(z)}{(n-1)
q^{(1)}(z)}\cdot\frac{q^{(3)}_n(z)}{(n-2)q^{(2)}_{n}(z)}\cdots\nonumber\\
&{}&
\cdots\frac{q^{(j-1)}_n(z)}{(n-j+2)q^{(j-2)}_n(z)}\cdot
\frac{q^{(j)}_n(z)}{(n-j+1)q^{(j-1)}_n(z)}\nonumber\\ & = &
\frac{q^{(j)}_n(z)}{n(n-1)\cdots (n-j+1)q_n(z)},\nonumber
\end{eqnarray}
and thus our basic assumption implies
\begin{eqnarray}\label{equationClimit}
C^j(z)=\lim_{n\to\infty}\prod_{i=0}^{j-1}C_{n,i}(z)=
\lim_{n\to\infty}\frac{q^{(j)}_n(z)}{n(n-1)\cdots (n-j+1)q_n(z)}.
\end{eqnarray}
 
In the above notation consider the eigenvalue equation $Tp_n(z)=\lambda_np_n(z)$, where
 the eigenvalue $\lambda_n$ is given by
\begin{displaymath}
\lambda_n=\sum_{j=1}^{k}\alpha_{j,j}\frac{n!}{(n-j)!}=
\sum_{j=1}^{j_0}\alpha_{j,j}\frac{n!}{(n-j)!}=
\sum_{j=1}^{j_0}\alpha_{j,j}n(n-1)\cdots(n-j+1).
\end{displaymath}
Clearly this sum ends at $j_0$ since $\alpha_{j,j}=0$ for all $j>j_0$ 
by definition of $j_0$ as the largest $j$ for which $\deg Q_j=j$.
 We then have
\begin{displaymath}
Tp_n(z)=\lambda_np_n(z)
\end{displaymath}
\begin{displaymath}
\Leftrightarrow
\end{displaymath}
\begin{displaymath}
\sum_{j=1}^{k}\bigg(\sum_{i=0}^{\deg Q_j}\alpha_{j,i}z^{i}\bigg)p^{(j)}_n(z)
=\sum_{j=1}^{j_0}\alpha_{j,j}n(n-1)\cdots (n-j+1)p_n(z).
\end{displaymath}
Substituting $z=n^dz$ in this equation we obtain
\begin{equation*}\label{ekvation}
\sum_{j=1}^{k}\bigg(\sum_{i=0}^{\deg Q_j}\alpha_{j,i}n^{di}z^{i}\bigg)p^{(j)}_n(n^dz)
=\sum_{j=1}^{j_0}\alpha_{j,j}n(n-1)\cdots (n-j+1)p_n(n^dz),
\end{equation*}
and with $q_n(z)=p_n(n^dz)$ we get
\begin{displaymath}
\sum_{j=1}^{k}\bigg(\sum_{i=0}^{\deg Q_j}\alpha_{j,i}\frac{z^{i}}{n^{d(j-i)}}\bigg)q^{(j)}_n(z)
=\sum_{j=1}^{j_0}\alpha_{j,j}n(n-1)\cdots (n-j+1)q_n(z).
\end{displaymath}
Dividing this equation by $\frac{n!}{(n-j_0)!}q_n(z)=n(n-1)\cdots (n-j_0+1)q_n(z)$ we get 
\begin{eqnarray}\label{scaledpoleigenvalueproblem}
 & &\sum_{j=1}^{k}\bigg(\sum_{i=0}^{\deg Q_j}\alpha_{j,i}\frac{z^{i}}{n^{d(j-i)}}\bigg)\frac{q^{(j)}_n(z)}{n(n-1)\cdots (n-j_0+1)q_n(z)}=\nonumber\\
&=&\sum_{j=1}^{j_0}\alpha_{j,j}\frac{n(n-1)\cdots (n-j+1)}{n(n-1)\cdots (n-j_0+1)}.
\end{eqnarray}
Consider the right-hand side of (\ref{scaledpoleigenvalueproblem}). Since $j\leq j_0$ all terms for which $j<j_0$ (if not already zero, which is the  case if $\alpha_{j,j}=0$, i.e. if $\deg Q_j<j$) tend to zero when $n\to\infty$, and therefore the limit of the right-hand side of (\ref{scaledpoleigenvalueproblem}) equals 
\begin{displaymath}
\lim_{n\to\infty}\sum_{j=1}^{j_0}\alpha_{j,j}\frac{n(n-1)\cdots (n-j+1)}{n(n-1)\cdots (n-j_0+1)}
=\alpha_{j_0,j_0}=1,
\end{displaymath}
since we assumed that $Q_{j_0}$ is monic. 
Now consider the $j$th term in the sum on the left-hand side of (\ref{scaledpoleigenvalueproblem}). It equals 
\begin{eqnarray*}
&{}&\sum_{i=0}^{\deg Q_j}\alpha_{j,i}\frac{z^{i}}{n^{d(j-i)}}\cdot\frac{q^{(j)}_n(z)}{n(n-1)\cdots (n-j_0+1)q_n(z)}=\\
&=&\sum_{i=0}^{\deg Q_j}\alpha_{j,i}\frac{z^{i}}{n^{d(j-i)}}\cdot\frac{q^{(j)}_n(z)}
{n(n-1)\cdots (n-j+1)q_n(z)}\cdot\frac{n(n-1)\cdots (n-j+1)}{n(n-1)\cdots (n-j_0+1)}\\
&=&\sum_{i=0}^{\deg Q_j}\alpha_{j,i}\frac{z^{i}}{n^{d(j-i)}}\cdot 
\prod_{i=0}^{j-1}C_{n,i}(z)
\cdot\frac{n(n-1)\cdots (n-j+1)}{n(n-1)\cdots (n-j_0+1)}\\
&=& \sum_{i=0}^{\deg Q_j}\alpha_{j,i}\frac{z^{i}}{n^{d(j-i)+j_0-j}}\cdot 
\prod_{i=0}^{j-1}C_{n,i}(z)
\cdot\frac{n(n-1)\cdots (n-j+1)}{n^j}\frac{n^{j_0}}{n(n-1)\cdots (n-j_0+1)}.\\
\end{eqnarray*}
Taking the limit and using the basic assumption (\ref{equationClimit}) we obtain 
\begin{eqnarray*}
&{}& \lim_{n\to\infty}\sum_{i=0}^{\deg Q_j}\alpha_{j,i}\frac{z^{i}}{n^{d(j-i)}}\cdot\frac{q^{(j)}_n(z)}{n(n-1)\cdots (n-j_0+1)q_n(z)}\\
&=&\lim_{n\to\infty}\sum_{i=0}^{\deg Q_j}\alpha_{j,i}
\frac{z^{i}}{n^{d(j-i)+j_0-j}}C^j(z) 
\end{eqnarray*}
for the $j$th term and 
thus, taking the limit of the left-hand side of (\ref{scaledpoleigenvalueproblem}) we get 
\begin{eqnarray*}
& &\lim_{n\to\infty}\sum_{j=1}^{k}\bigg(\sum_{i=0}^{\deg Q_j}\alpha_{j,i}\frac{z^{i}}{n^{d(j-i)}}\bigg)\frac{q^{(j)}_n(z)}{n(n-1)\cdots (n-j_0+1)q_n(z)}\\
&=&\lim_{n\to\infty}\sum_{j=1}^{k}\bigg(\sum_{i=0}^{\deg Q_j}\alpha_{j,i}
\frac{z^{i}}{n^{d(j-i)+j_0-j}}\bigg)C^j(z).
\end{eqnarray*}
Adding up, the following equation is satisfied by 
$C(z)$ for $z$'s away from the support of $\mu$:
\begin{eqnarray}\label{preliminaryCekva}
\lim_{n\to\infty}\sum_{j=1}^{k}\bigg(\sum_{i=0}^{\deg Q_j}\alpha_{j,i}
\frac{z^{i}}{n^{d(j-i)+j_0-j}}\bigg)C^j(z)=1.
\end{eqnarray}
In order to make (\ref{preliminaryCekva}) a well-defined algebraic equation, i.e. to avoid infinities in the denominator when $n\to\infty$, we must impose the following condition on the real number $d$ in the exponent of $n$, namely  
\begin{displaymath}
d(j-i)+j_0-j\geq 0\quad\Leftrightarrow\quad d\geq\frac{j-j_0}{j-i}
\end{displaymath}
 for all $j\in[1,k]$ and all $i\in[0,\deg Q_j]$. Therefore we take $d=
\max_{j\in[1,k]\atop i\in[0,\deg Q_j]}\big(\frac{j-j_0}{j-i}\big)$, but this maximum is clearly obtained for the maximal value of $i$ for any given $j$, so we may as well put $i=\deg Q_j$. Our condition then becomes $d=\max_{j\in[1,k]}\big(\frac{j-j_0}{j-\deg Q_j}\big)$, and clearly the maximum is taken only over 
$j$ for which $Q_j(z)$ is not identically zero.  
Finally we observe that since $T$ is degenerate we have $j_0<k$ and thus we need only take this maximum over $j\in[j_0+1,k]$, since 
there always exists a positive value on $d$ for any operator of the type we consider. Thus our condition becomes:
\begin{displaymath}
d=\max_{j\in[j_0+1,k]}\bigg(\frac{j-j_0}{j-\deg Q_j}\bigg).
\end{displaymath}\\

\qquad\qquad\qquad\qquad \textbf{Corollary of Main Conjecture.}\\\\
In the above notation (recall that $\alpha_{j_0,j_0}=1$ by monicity of $Q_{j_0}$), 
the following well-defined algebraic equation follows immediately 
from inserting $d$ as defined in Main Conjecture into equation (\ref{preliminaryCekva}) and letting $n\to\infty$:\\\\
\textbf{Corollary.} \textit{The Cauchy transform $C(z)$ of the asymptotic root measure $\mu$ of the scaled eigenpolynomial $q_n(z)=p_n(n^dz)$ of an arbitrary exactly-solvable operator $T$ as above satisfies the following algebraic equation for almost all complex $z$ in the usual Lebesgue measure on $\mathbb{C}$}: 
\begin{displaymath}
z^{j_0}C^{j_0}(z)+ \sum_{j\in A}\alpha_{j,\deg Q_{j}}z^{\deg Q_j}
C^j(z)=1, 
\end{displaymath}
\textit{where $A$ is the set consisting of all $j$ for which the maximum\\
 $d:=\max_{j\in [j_0+1,k]}\big(\frac{j-j_0}{j-\deg Q_j}\big)$ is attained,
 i.e. 
$A=\{j:(j-j_0)/(j-\deg Q_j)=d\}$.}\\\\
For detailed calculations see Section 5.1.
\section{Numerical evidence}
\subsection{Evidence for Main Conjecture}
In the table on the last page of this section we present 
numerical evidence on the growth of $r_n=\max\{|z|:p_n(z)=0\}$ which supports the choice of $d$ in Main Conjecture.
We have performed similar computer experiments for a large number of other degenerate exactly-solvable operators, and the results are in all cases consistent with Main Conjecture. Next we present some typical pictures on the zero distribution of the appropriately scaled eigenpolynomial $q_n(z)=p_n(n^dz)$ 
for some degenerate exactly-solvable operators, where $p_n$ denotes the $n$th degree unique and monic polynomial eigenfunction of the given operator $T$. Conjecturally the zeros of $q_n(z)$ are contained in a compact set when $n\to\infty$.\\\\\\\\
\textit{Fig.1:} $T_1=zD+zD^2+zD^3+zD^4+zD^5$.\\
\begin{tabular}{ccccc}
\includegraphics[width=2.5cm]{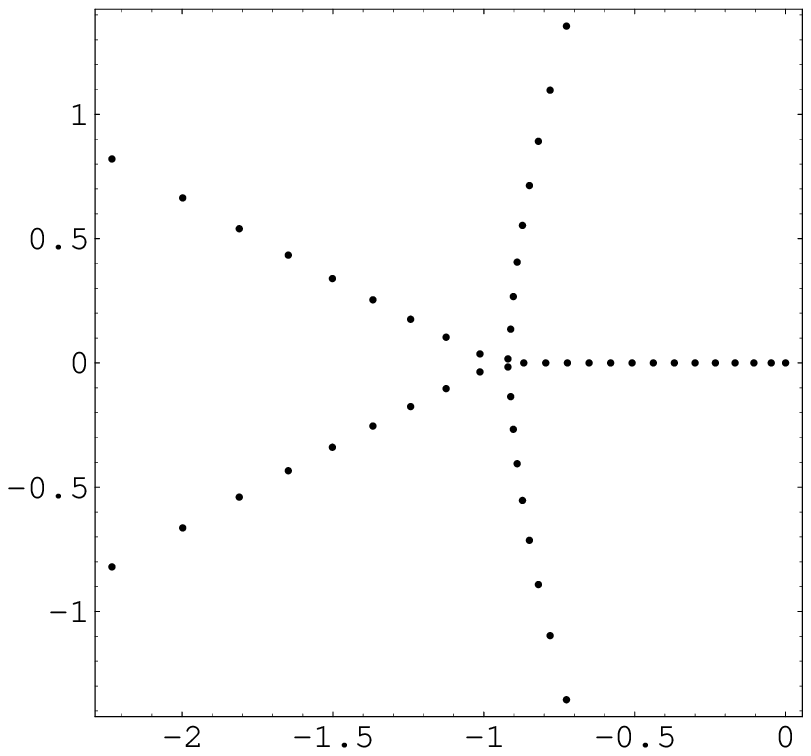} & &
\includegraphics[width=2.5cm]{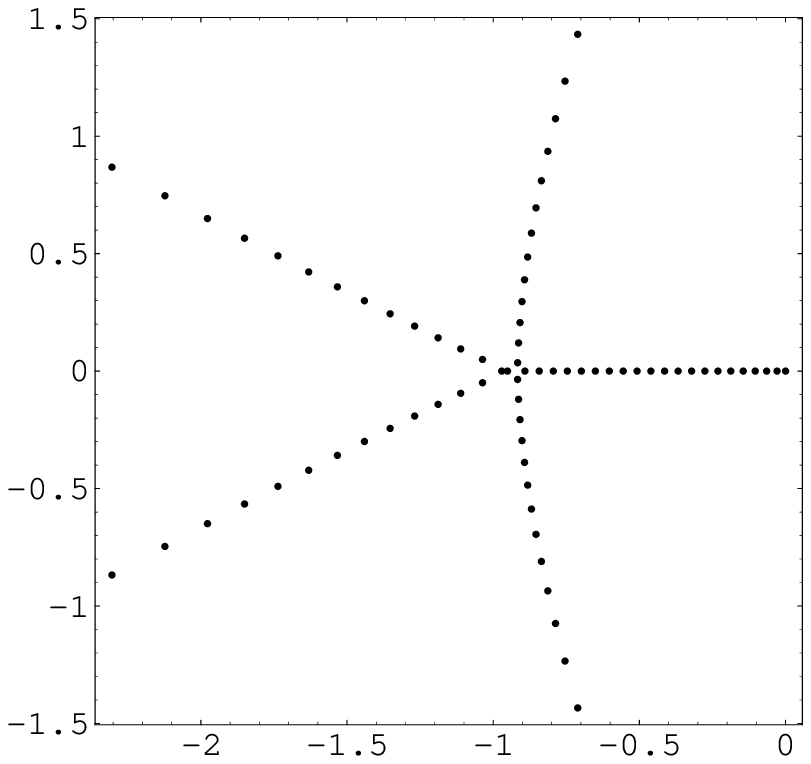} & &
\includegraphics[width=2.5cm]{T1bild100_gr1.eps}\\
roots of & & roots of & & roots of\\
 $q_{50}(z)=p_{50}(50z)$ & & $q_{75}(z)=p_{75}(75z)$
 & & $q_{100}(z)=p_{100}(100z)$ 
\end{tabular}\\
\\\\
\textit{Fig.2:} $T_2=z^2D^2+D^7$.\\
\begin{tabular}{ccccc}
\includegraphics[width=2.5cm]{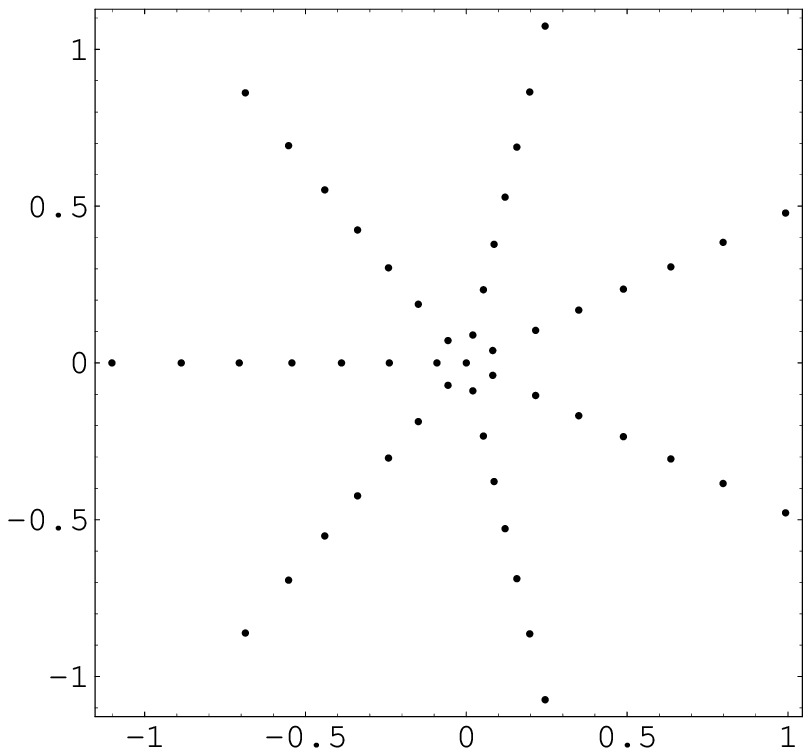} & &
\includegraphics[width=2.5cm]{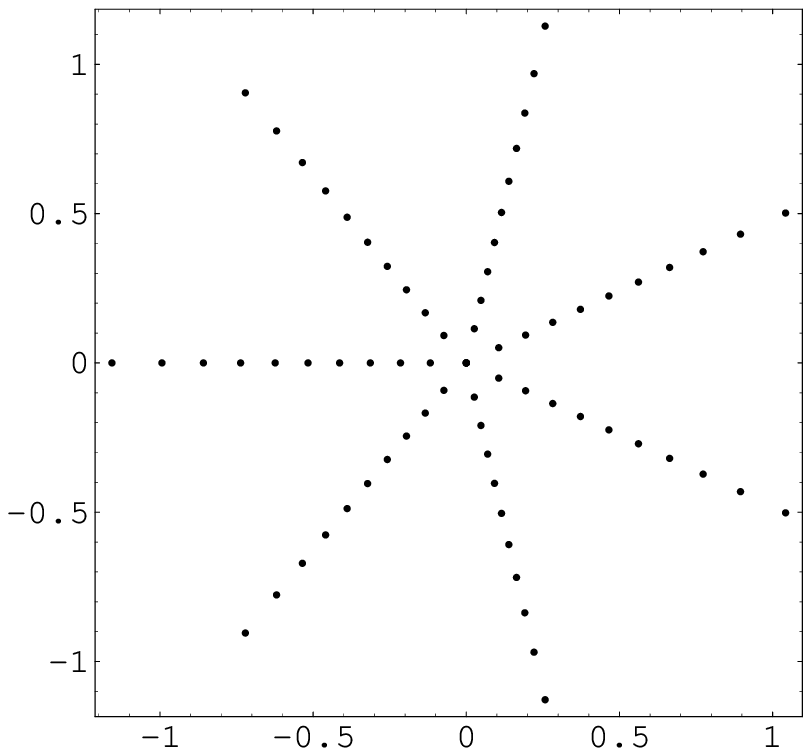} & &
\includegraphics[width=2.5cm]{T2bild100_gr1.eps}\\
roots of & & roots of & & roots of\\
 $q_{50}(z)=p_{50}(50^{5/7}z)$ & &
 $q_{75}(z)=p_{75}(75^{5/7}z)$
 & & $q_{100}(z)=p_{100}(100^{5/7}z)$
\end{tabular}\\
\\\\
\textit{Fig.3:} $T_3=z^3D^3+z^2D^4+zD^5$.\\
\begin{tabular}{ccccc}
\includegraphics[width=2.5cm]{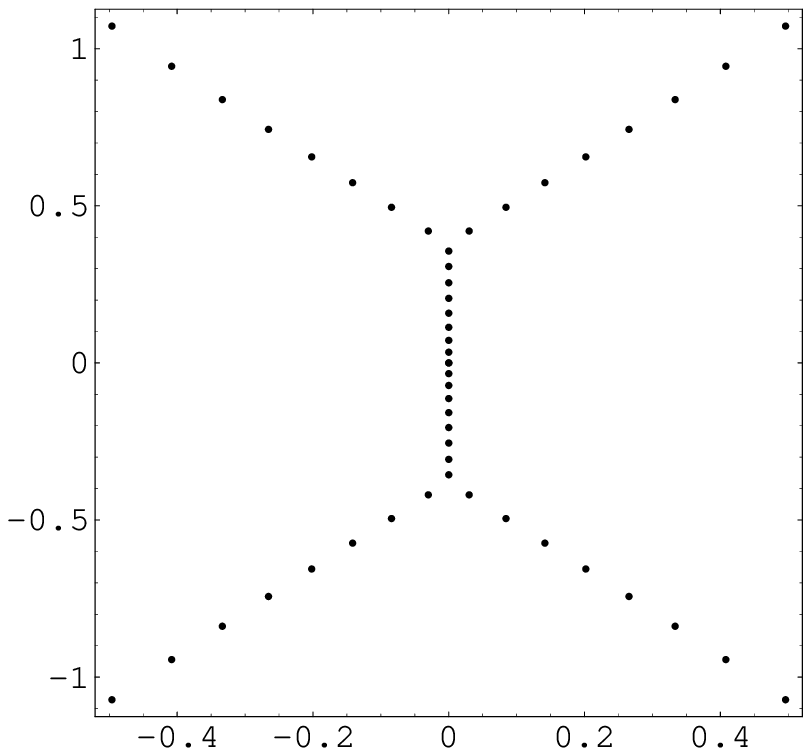} & &
\includegraphics[width=2.5cm]{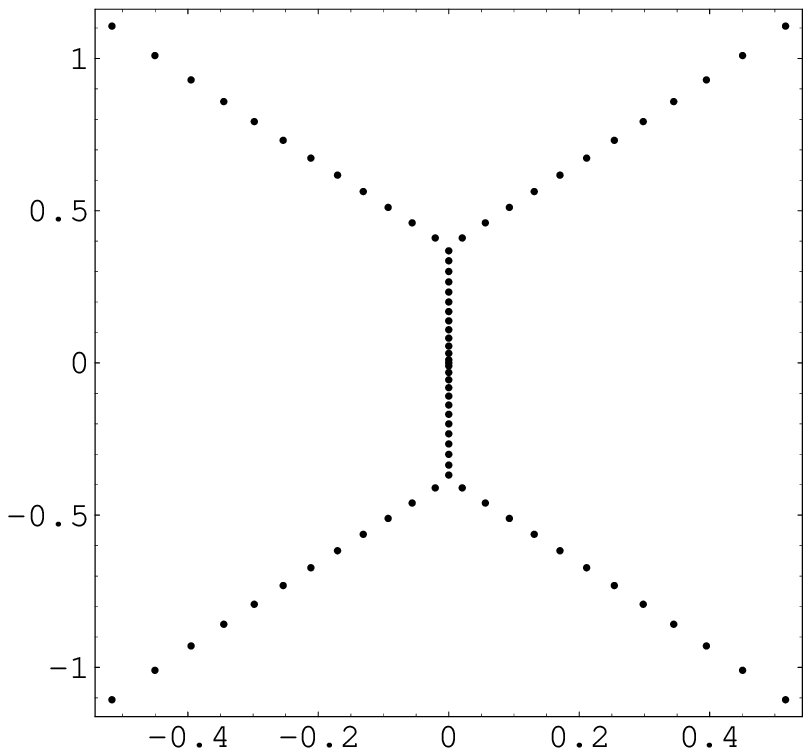} & &
\includegraphics[width=2.5cm]{T3bild100_gr1.eps}\\
roots of & & roots of & & roots of\\
 $q_{50}(z)=p_{50}(50^{1/2}z)$ & &
 $q_{75}(z)=p_{75}(75^{1/2}z)$ & & $q_{100}(z)=p_{100}(100^{1/2}z)$
\end{tabular}
\subsection{On the basic assumption}
Let us show some examples of pictures supporting the basic assumption, namely that the 
root measures of $q_n(z)$ and its derivatives up to the $k$th
 order exist when $n\to\infty$ and are all weakly convergent to the same measure $\mu$.\\\\
\textit{Fig. 1:} $T_7=zD+D^3$ and $q_n(z)=p_n(n^{2/3}z)$.\\ 
\begin{tabular}{cc}
\includegraphics[width=2.5cm]{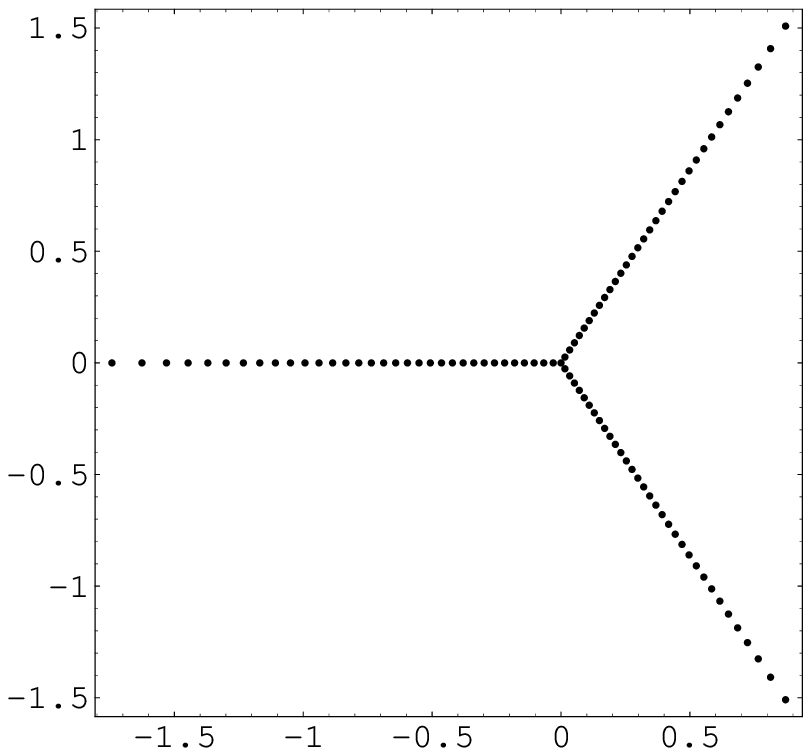} & 
\includegraphics[width=2.5cm]{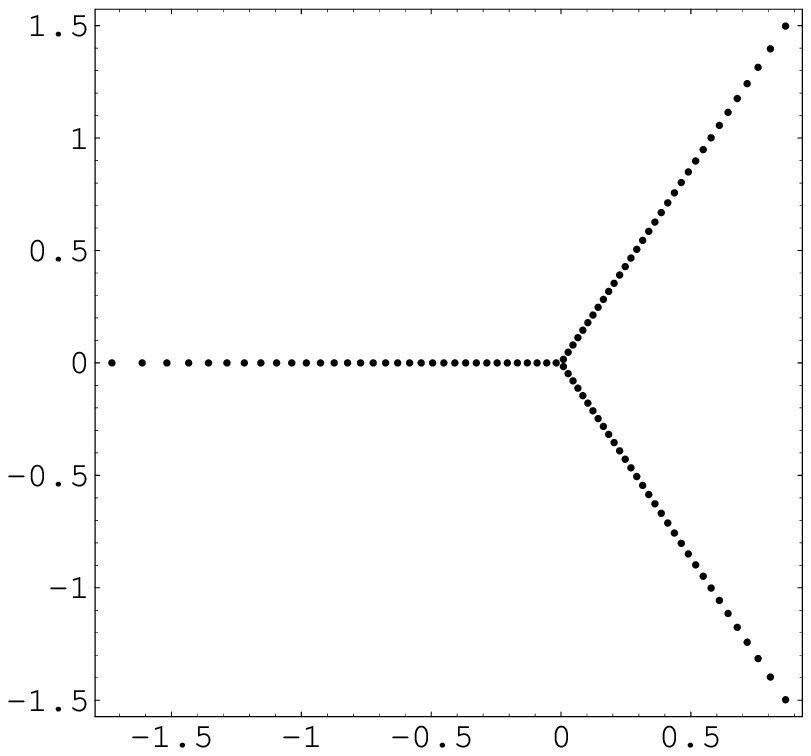}\\
roots of $q_{100}(z)$ & 
roots of $q^{\prime}_{100}(z)$
\end{tabular}
\begin{tabular}{cc}
\includegraphics[width=2.5cm]{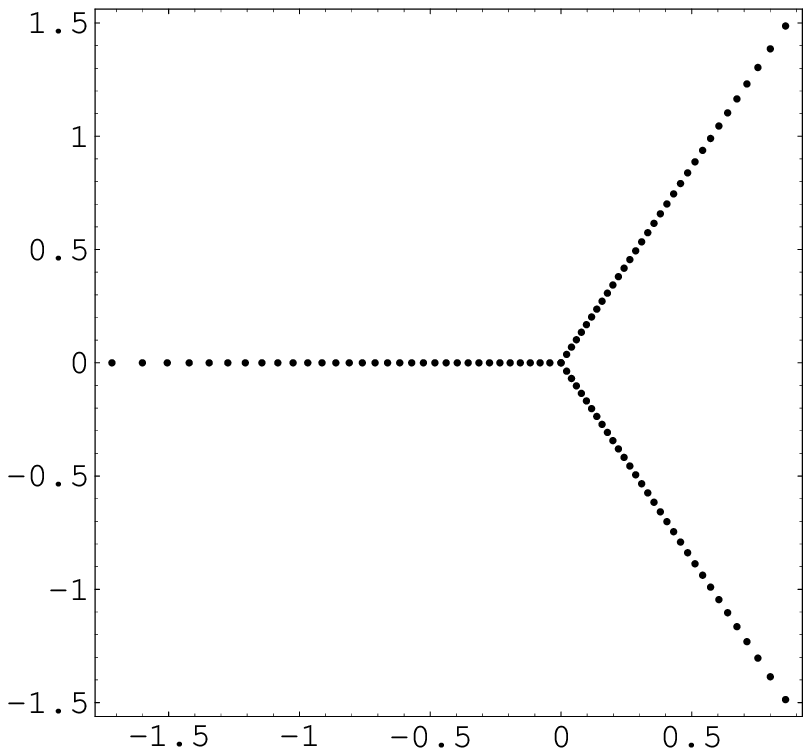} & 
\includegraphics[width=2.5cm]{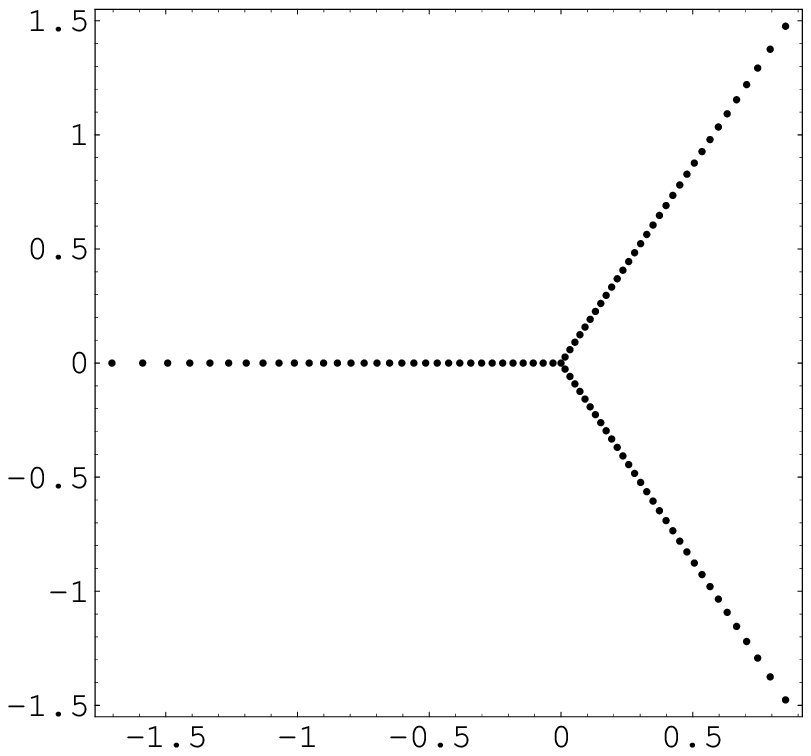}\\
roots of $q^{\prime\prime}_{100}(z)$  & 
roots of $q^{\prime\prime\prime}_{100}(z)$
\end{tabular}\\\\\\
\textit{Fig. 2:} $T_{9}=zD+zD^4+z^3D^7$ and $q_n(z)=p_n(n^{3/2}z)$.\\
\begin{tabular}{cccc}
\includegraphics[width=2.5cm]{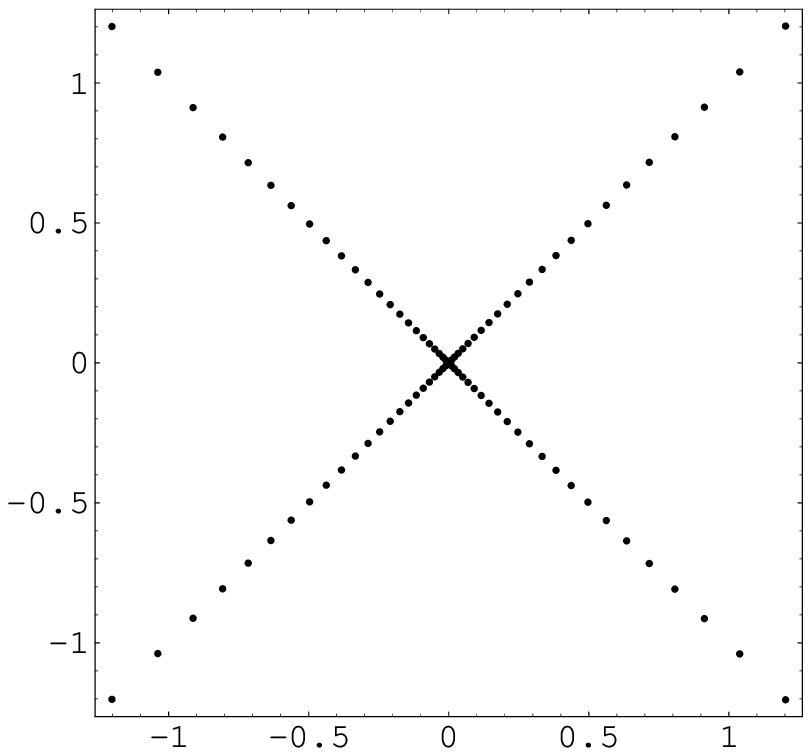} & 
\includegraphics[width=2.5cm]{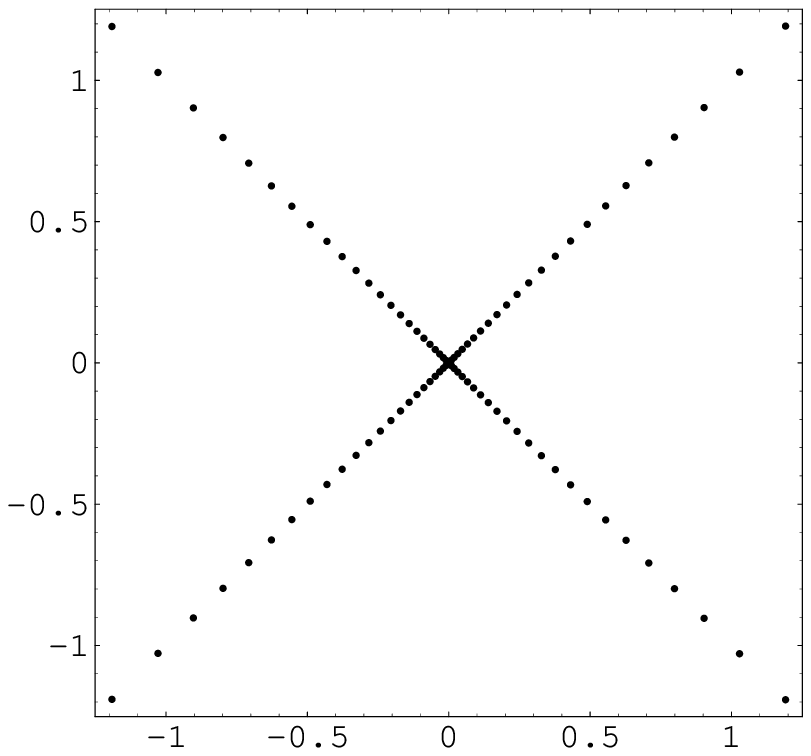} & 
\includegraphics[width=2.5cm]{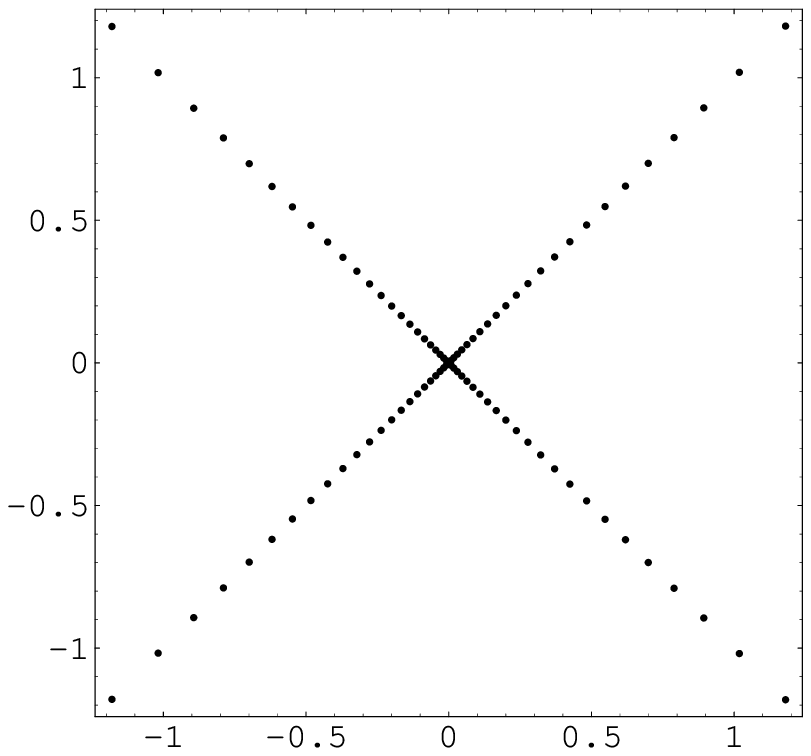} &
\includegraphics[width=2.5cm]{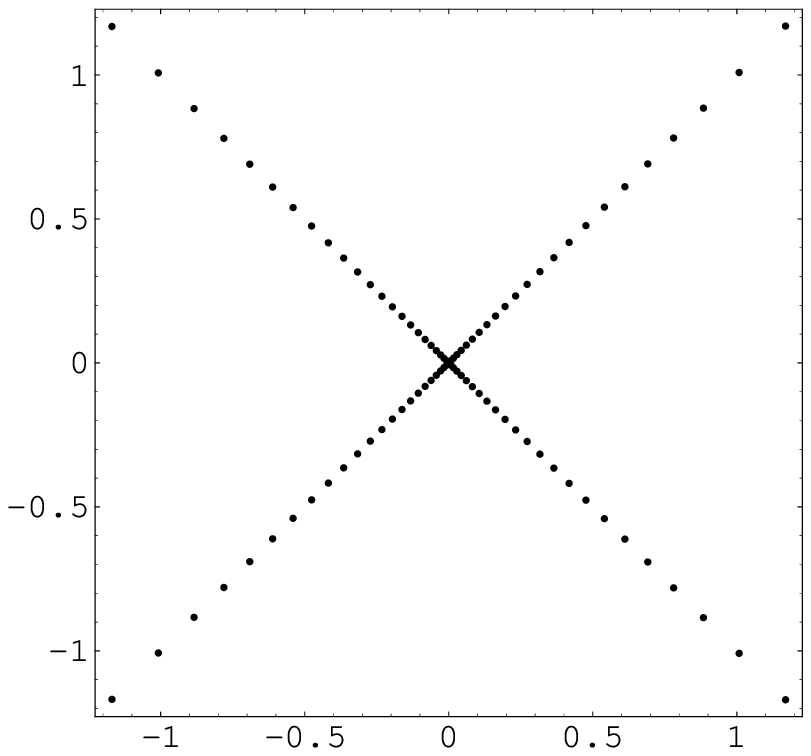}\\
roots of $q_{100}(z)$  & 
roots of $q^{\prime}_{100}(z)$ & roots of $q^{\prime\prime}_{100}(z)$ &
roots of  $q^{\prime\prime\prime}_{100}(z)$
\end{tabular}\\\\
\begin{tabular}{cccc}
\includegraphics[width=2.5cm]{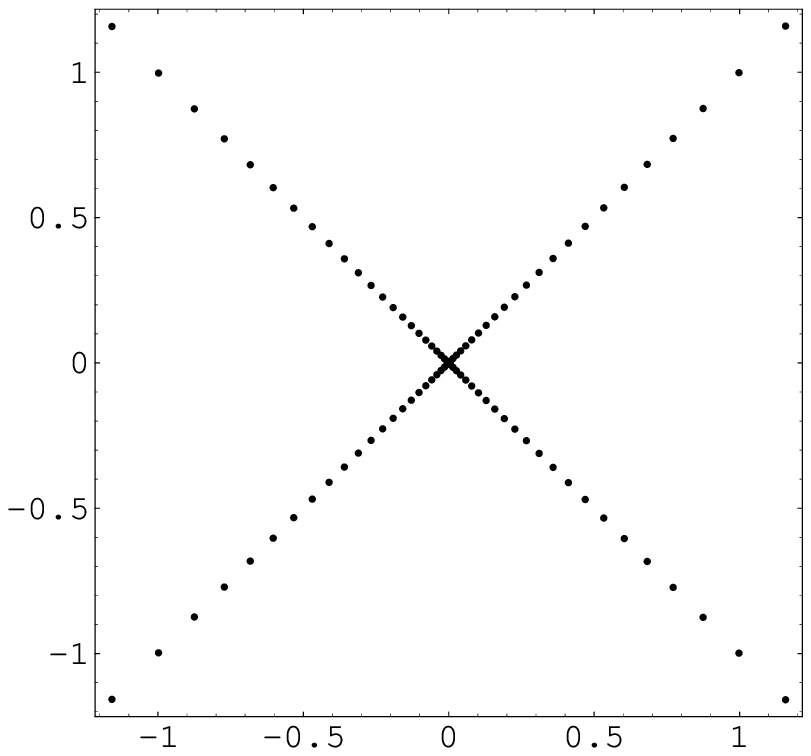} & 
\includegraphics[width=2.5cm]{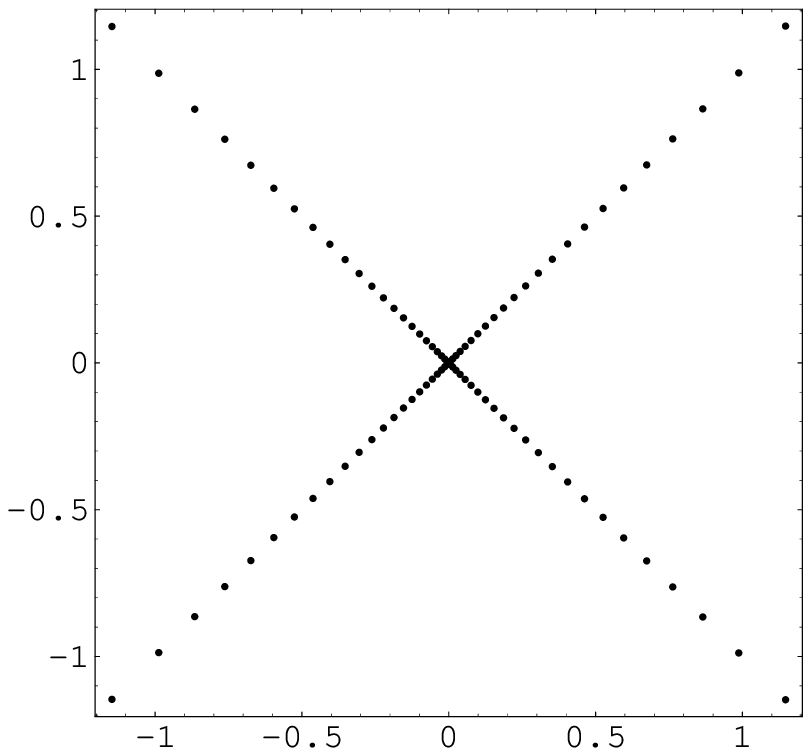} & 
\includegraphics[width=2.5cm]{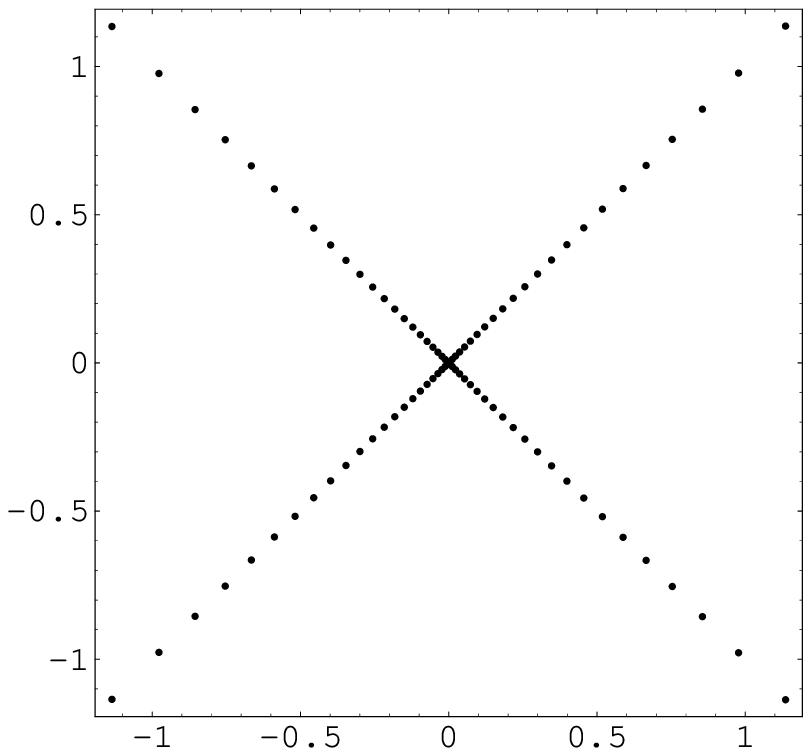} &
\includegraphics[width=2.5cm]{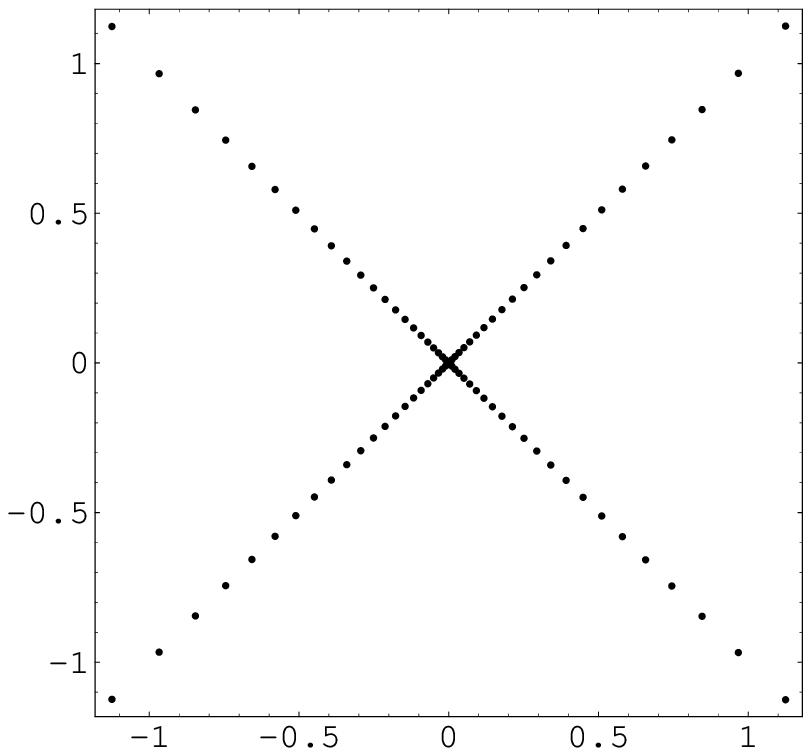}\\
roots of $q^{(\rm{iv})}_{100}(z)$  & 
roots of $q^{(\rm{v})}_{100}(z)$ & roots of $q^{(\rm{vi})}_{100}(z)$ &
roots of  $q^{(\rm{vii})}_{100}(z)$
\end{tabular}\\
\subsection{On the Corollary of Main Conjecture}
 The algebraic equation in the corollary of Main Conjecture satisfied by the Cauchy transform 
of the asymptotic root measure of the scaled eigenpolynomial  
 indicates that the asymptotic zero distribution depends only on the 
term $z^{j_0}D^{j_0}$ and the term(s) $\alpha_{j,\deg Q_j}z^{\deg Q_j}D^j$ for which $d=\max_{j\in[j_0+1,k]}\big(\frac{j-j_0}{j-\deg Q_j}\big)$ is attained.\footnote{Recall that we normalized $T$ by letting $Q_{j_0}$ be monic, i.e. 
$\alpha_{j_0,j_0}=1$.} Thus any term $Q_jD^j$ in $T$ for which $j<j_0$ or such that $(j-j_0)/(j-\deg Q_j)<d$ is (conjecturally) irrelevant for the zero distribution when $n\to\infty$.\\

To illustrate this fact we now present some pictures of the zero distributions of the scaled eigenpolynomials for some \textit{distinct} operators for which the Cauchy transforms $C$ of the corresponding scaled eigenpolynomials $q_n$ satisfy the \textit{same} equation when $n\to\infty$, namely the equation in the Corollary of Main Conjecture in   Section 3.\\

As a first example consider the operator $T_4=z^3D^3+z^2D^5$. Here $d=
\max_{j\in[j_0+1,k]}\big(\frac{j-j_0}{j-\deg Q_j}\big)=(5-3)/(5-2)=2/3$, the corresponding scaled eigenpolynomial is $q_n(z)=
p_n(n^{2/3}z)$, and we have  
$z^3C^3+z^2C^5=1$ for the Cauchy transform of $q_n$ when $n\to\infty$. Now consider the slightly modified operator $\widetilde{T}_4=z^2D^2+z^3D^3+zD^4+z^2D^5+D^6$ and note that  $d$ is obtained again (only) for $j=5$ (for $j=4$ we have $(4-3)/(4-1)=1/3<2/3$ and for $j=6$ we have $(6-3)/(6-0)=3/6=1/2<2/3$). We therefore obtain the same 
Cauchy transform equation as for $T_4$, and hence the terms $z^2D^2$, $zD^4$ and $D^6$ in $\widetilde{T}_4$ can be considered as irrelevant for the zero distribution for sufficiently large $n$. The pictures below clearly illustrate this.\\\\
\begin{tabular}{ccc}
\includegraphics[width=2.6cm]{T4bild100_gr1.eps} &  &
\includegraphics[width=2.6cm]{T4bildcompare100_gr1.eps}\\
$T_4=z^3D^3+z^2D^5$, & & $\widetilde{T}_4=z^2D^2+z^3D^3+zD^4+z^2D^5+D^6$,\\
roots of $q_{100}(z)=p_{100}(100^{2/3}z)$ & & roots of $q_{100}(z)=p_{100}(100^{2/3}z)$
\end{tabular}\\\\

However, instead of $D^6$, we may add the more ''disturbing'' term $zD^6$ to $T_4$. Note that for the operator $\overline{T}_4=z^2D^2+z^3D^3+zD^4+z^2D^5+zD^6$ for $j=6$ we have $(6-3)/(6-1)=3/5=0.6<2/3$. 
Adding any term $Q_jD^j$ such that $(j-j_0)/(j-\deg Q_j)<d$ to a given operator, it is clear that the 
closer the value of $(j-j_0)/(j-\deg Q_j)$ is to $d$ (in this case $2/3$), the more disturbing it is in the sense that it requires larger $n$ for the corresponding zero distributions to coincide. See pictures below.\\\\
\begin{tabular}{ccc}
\includegraphics[width=2.6cm]{T4bild100_gr1.eps} &
\includegraphics[width=2.6cm]{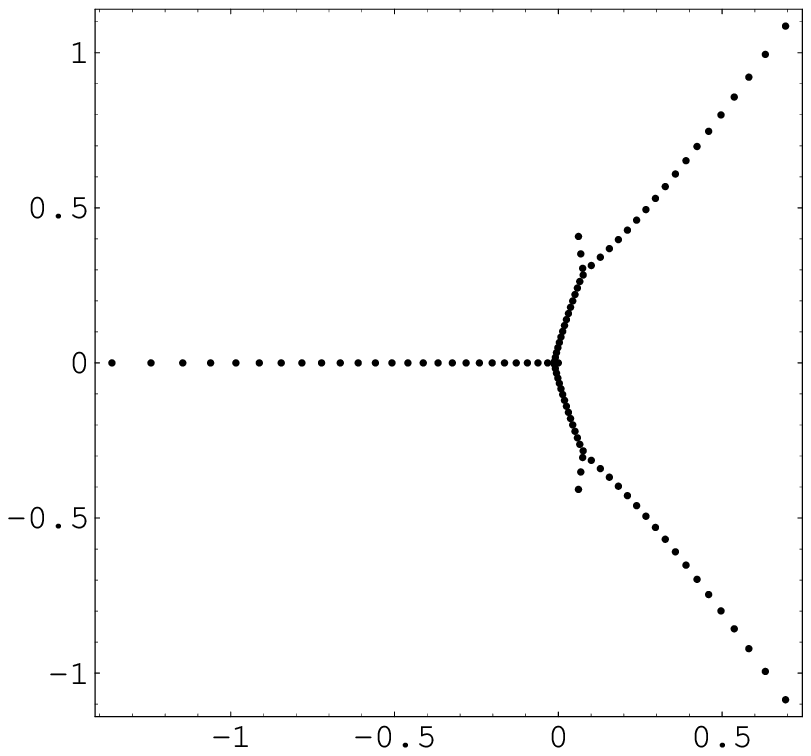}
& \includegraphics[width=2.6cm]{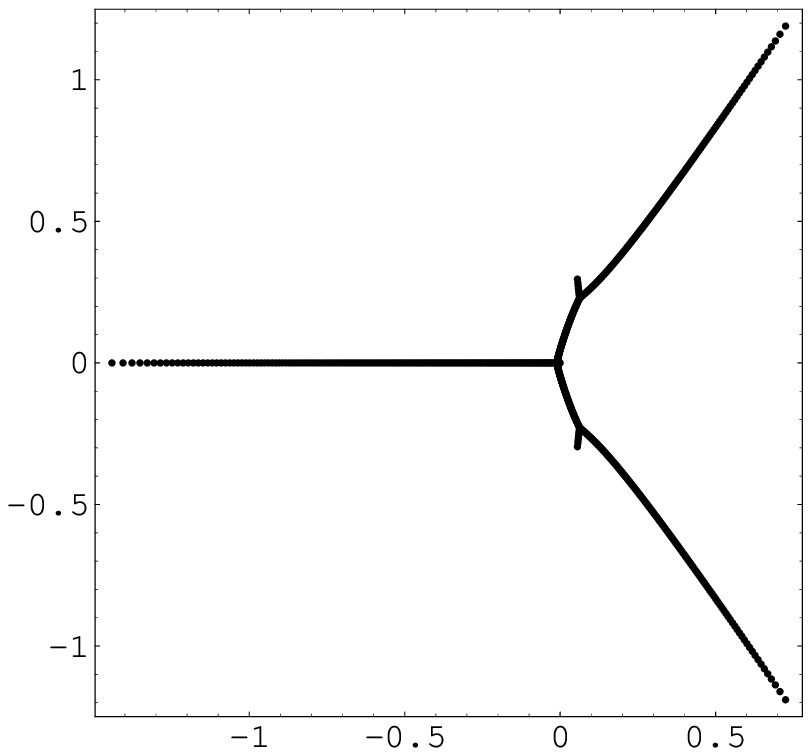}
\\
$T_4$, roots of & $\overline{T}_4$, roots of  & $\overline{T}_4$, roots of\\ 
 $q_{100}(z)=p_{100}(100^{2/3}z)$ & $q_{100}(z)=p_{100}(100^{2/3}z)$ & $q_{600}(z)=p_{600}(600^{2/3}z)$
\end{tabular}\\\\\\
Increasing $n$ however, experiments indicate that the zero distributions 
of the scaled eigenpolynomials of $T_4$ and $\overline{T}_4$ coincide, as they (conjecturally) should.\\

 As a second example, consider the operators $T_5=z^5D^5+z^4D^6+z^2D^8$ and 
$\widetilde{T}_5=z^2D^2+z^5D^5+z^4D^6+zD^7+z^2D^8$, whose scaled eigenpolynomials $q_n(z)=p_n(n^{1/2}z)$ both satisfy the Cauchy transform equation $z^5C^5+z^4C^6+z^2C^8=1$ when $n\to\infty$. In the pictures below we see that the terms $z^2D^2$ and $zD^7$ of $\widetilde{T}_5$ seem to have no effect on the zero distribution for large $n$.\\\\ 
\begin{tabular}{cc}
\includegraphics[width=2.6cm]{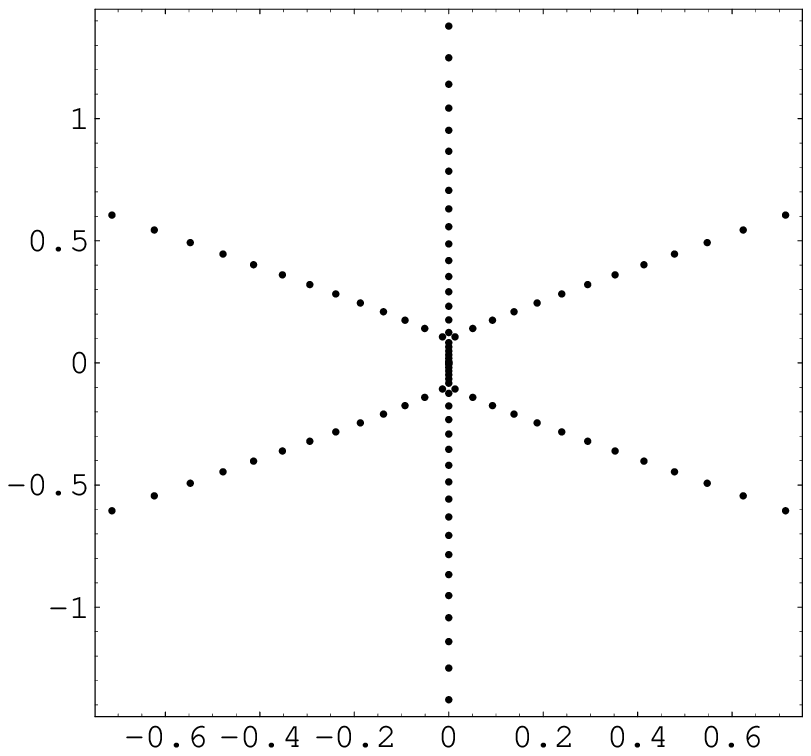} &
\includegraphics[width=2.6cm]{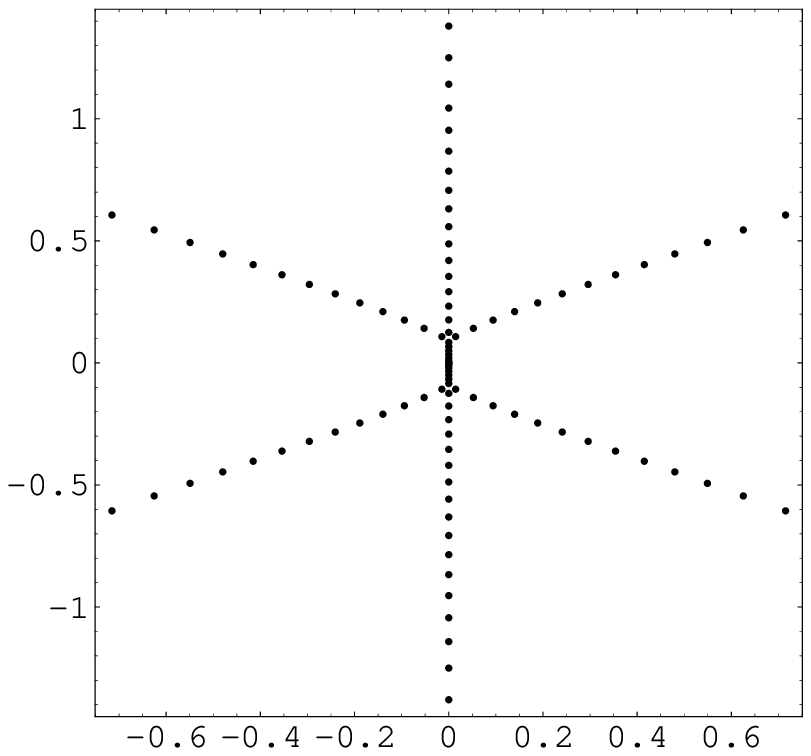}\\
$T_5=z^5D^5+z^4D^6+z^2D^8$, & $\widetilde{T}_5=z^2D^2+z^5D^5+z^4D^6+zD^7+z^2D^8$,\\
roots of $q_{100}(z)=p_{100}(100^{1/2}z)$ & roots of $q_{100}(z)=p_{100}(100^{1/2}z)$
\end {tabular}\\\\

 Finally note that for $j_0$ and for any $j$ for which $d$ is attained, it is only the
highest degree term $\alpha_{j,\deg Q_j}z^{\deg Q_j}$ of $Q_j$ that is
involved in the Cauchy transform equation. 
 Consider for example the following case, where adding lower degree terms 
to $\alpha_{j,\deg Q_j}z^{\deg Q_j}$
in the (relevant) $Q_j$ seems to have no effect on the zero distribution for large $n$. Below 
$T_6=z^3D^3+z^2D^6$, and 
$\widetilde{T}_6=[(1+13i)+(24i-3)z+11iz^2+z^3]D^3+[(22i-13)+(-9-14i)z+z^2]D^6$.
Note the difference in scaling between the pictures.\\\\
\begin{tabular}{ccc}
\includegraphics[width=2.6cm]{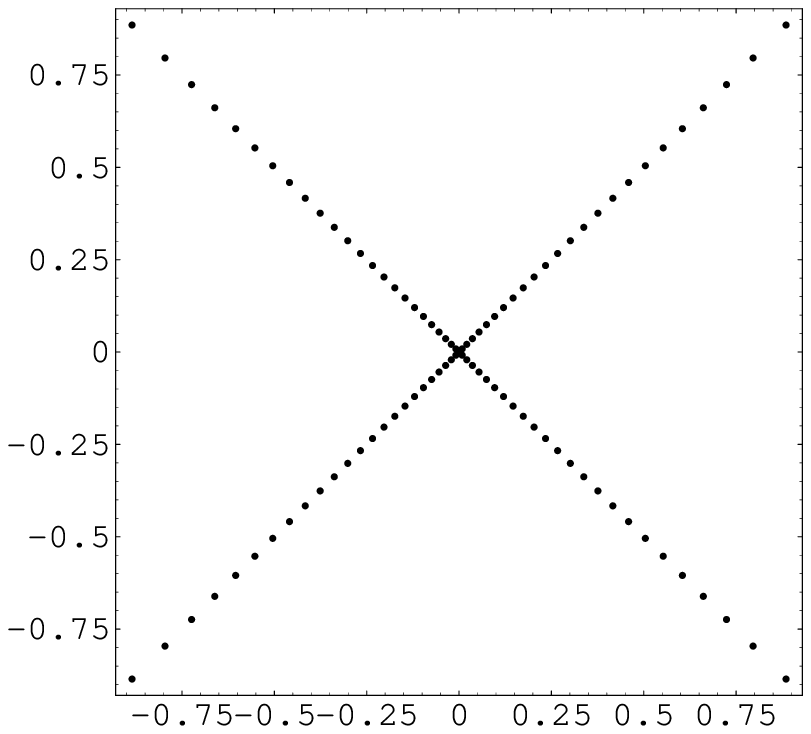} &
\includegraphics[width=2.6cm]{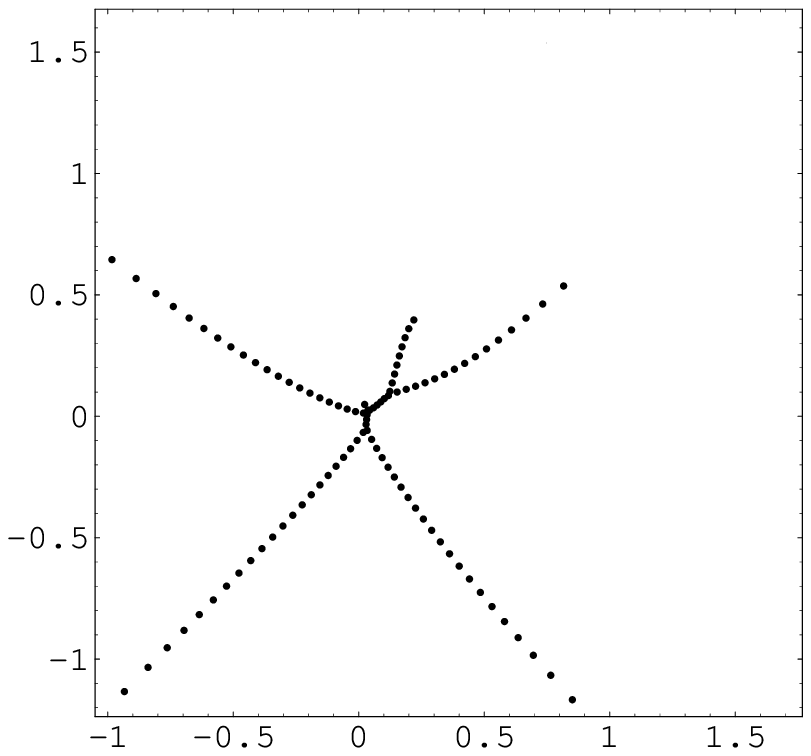} &
\includegraphics[width=2.6cm]{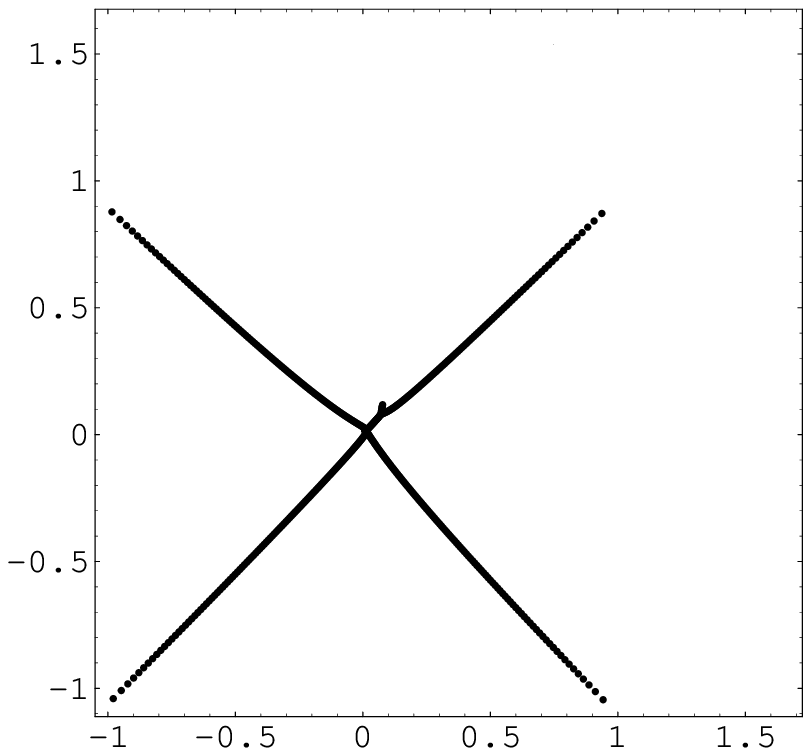}\\
$T_6$, roots of  & $\widetilde{T}_6$, roots of & 
$\widetilde{T}_6$, roots of\\
 $q_{100}(z)=p_{100}(100^{3/4}z)$ & $q_{100}(z)=p_{100}(100^{3/4}z)$ 
& $q_{500}(z)=p_{500}(500^{3/4}z)$
\end{tabular}\\

\label{tabell}
\begin{tabular}{|l|l|c|l|}
\hline
Operator & $n$ & $r_n$ experimental & 
$r_n$ conjectured\\
\hline
\hline
 & 50 & $2.7\cdot 50^{0.967595}$ & $c_{1}\cdot 50^{1}$\\
\cline{2-4}
$T_1=zD+zD^2+zD^3+zD^4+zD^5$  & 100 & $2.7\cdot 100^{0.984180}$ & $c_{1}\cdot 100^{1}$\\
\cline{2-4}
 & 200 & $2.7\cdot 200^{0.992557}$ & $c_{1}\cdot 200^{1}$\\
\cline{2-4}
 & 250 & $2.7\cdot 250^{0.994272}$ & $c_{1}\cdot 250^{1}$\\
\hline\hline
 & 50 & $1.3\cdot 50^{0.671977}$ & $c_{2}\cdot 50^{5/7}$\\
\cline{2-4}
 & 100 & $1.3\cdot 100^{0.694847}$ & $c_{2}\cdot 100^{5/7}$\\
\cline{2-4}
$T_2=z^2D^2+D^7$ & 200 & $1.3\cdot 200^{0.706226}$ & $c_{2}\cdot 200^{5/7}$\\
\cline{2-4}
 & 300 & $1.3\cdot 300^{0.710085}$ & $c_{2}\cdot 300^{5/7}$\\
\cline{2-4}
 & 400 & $1.3\cdot 400^{0.712043}$ & $c_{2}\cdot 400^{5/7}$\\
\hline\hline
 & 50 & $4/3\cdot 50^{0.469007}$ & $c_{3}\cdot 50^{1/2}$\\
\cline{2-4}
 & 100 & $4/3\cdot 100^{0.484824}$ & $c_{3}\cdot 100^{1/2}$\\
\cline{2-4}
$T_3=z^3D^3+z^2D^4+zD^5$ & 200 & $4/3\cdot 200^{0.492832}$ & $c_{3}\cdot 200^{1/2}$\\
\cline{2-4}
 & 300 & $4/3\cdot 300^{0.495592}$ & $c_{3}\cdot 300^{1/2}$\\
\cline{2-4}
 & 400 & $4/3\cdot 400^{0.497009}$ & $c_{3}\cdot 400^{1/2}$\\
\hline\hline
 & 50 & $1.4\cdot 50^{0.633226}$ & $c_{4}\cdot 50^{2/3}$\\
\cline{2-4}
 & 100 & $1.4\cdot 100^{0.652141}$ & $c_{4}\cdot 100^{2/3}$\\
\cline{2-4}
$T_4=z^3D^3+z^2D^5$ & 200 & $1.4\cdot 200^{0.661412}$ & $c_{4}\cdot 200^{2/3}$\\
\cline{2-4}
 & 300 & $1.4\cdot 300^{0.664511}$ & $c_{4}\cdot 300^{2/3}$\\
\cline{2-4}
 & 400 & $1.4\cdot 400^{0.666066}$ & $c_{4}\cdot 400^{2/3}$\\
\hline\hline
 & 50 & $1.4\cdot 50^{0.632811}$ & $\tilde{c}_{4}\cdot 50^{2/3}$\\
\cline{2-4}
 & 100 & $1.4\cdot 100^{0.651960}$ & $\tilde{c}_{4}\cdot 100^{2/3}$\\
\cline{2-4}
$\widetilde{T}_4=z^2D^2+z^3D^3+zD^4+z^2D^5+D^6$ & 200 & $1.4\cdot 200^{0.661332}$ & $\tilde{c}_{4}\cdot 200^{2/3}$\\
\cline{2-4}
 & 300 & $1.4\cdot 300^{0.664461}$ & $\tilde{c}_{4}\cdot 300^{2/3}$\\
\cline{2-4}
 & 400 & $1.4\cdot 400^{0.666030}$ & $\tilde{c}_{4}\cdot 400^{2/3}$\\
\hline\hline
 & 50 & $1.5\cdot 50^{0.462995}$ & $c_{5}
\cdot 50^{1/2}$\\
\cline{2-4}
 & 100 & $1.5\cdot 100^{0.481684}$ & $c_{5}\cdot 100^{1/2}$\\
\cline{2-4}
$T_5=z^5D^5+z^4D^6+z^2D^8$ & 200 & $1.5\cdot 200^{0.491066}$ & $c_{5}\cdot 200^{1/2}$\\
\cline{2-4}
 & 300 & $1.5\cdot 300^{0.494304}$ & $c_{5}\cdot 300^{1/2}$\\
\cline{2-4}
 & 400 & $1.5\cdot 400^{0.495971}$ & $c_{5}\cdot 400^{1/2}$\\
\hline\hline
 & 50 & $1.5\cdot 50^{0.463391}$ & $\tilde{c}_{5}
\cdot 50^{1/2}$\\
\cline{2-4}
 & 100 & $1.5\cdot 100^{0.481837}$ & $\tilde{c}_{5}\cdot 100^{1/2}$\\
\cline{2-4}
$\widetilde{T}_5=z^2D^2+z^5D^5+z^4D^6+zD^7+z^2D^8$ & 200 & $1.5\cdot 200^{0.491129}$ & 
$\tilde{c}_{5}\cdot 200^{1/2}$\\
\cline{2-4}
 & 300 & $1.5\cdot 300^{0.494342}$ & $\tilde{c}_{5}\cdot 300^{1/2}$\\
\cline{2-4}
 & 400 & $1.5\cdot 400^{0.495998}$ & $\tilde{c}_{5}\cdot 400^{1/2}$\\
\hline\hline
 & 50 & $1.4\cdot 50^{0.702117}$ & $c_{6}\cdot 50^{3/4}$\\
\cline{2-4}
 & 100 & $1.4\cdot 100^{0.725715}$ & $c_{6}\cdot 100^{3/4}$\\
\cline{2-4}
$T_6=z^3D^3+z^2D^6$ & 200 & $1.4\cdot 200^{0.737541}$ & $c_{6}\cdot 200^{3/4}$\\
\cline{2-4}
 & 300 & $1.4\cdot 300^{0.741614}$ & $c_{6}\cdot 300^{3/4}$\\
\cline{2-4}
 & 400 & $1.4\cdot 400^{0.743713}$ & $c_{6}\cdot 400^{3/4}$\\
\hline\hline
 & 50 & $1.4\cdot 50^{0.769260}$ & $\tilde{c}_{6}\cdot 50^{3/4}$\\
\cline{2-4}
$\widetilde{T}_6=[(1+13i)+(24i-3)z+11iz^2+z^3]D^3$ & 100 & $1.4\cdot 100^{0.760399}$ & 
$\tilde{c}_{6}\cdot 100^{3/4}$\\
\cline{2-4}$+[(22i-13)-(9+14i)z+z^2]D^6$
 & 200 & $1.4\cdot 200^{0.756161}$ & $\tilde{c}_{6}\cdot 200^{3/4}$\\
\cline{2-4}
 & 300 & $1.4\cdot 300^{0.754590}$ & $\tilde{c}_{6}\cdot 300^{3/4}$\\
\cline{2-4}
 & 400 & $1.4\cdot 400^{0.753765}$ & $\tilde{c}_{6}\cdot 400^{3/4}$\\
\hline
\end{tabular}\\
\section{Appendix} 
\subsection{Arriving at the Corollary of Main Conjecture} 
The algebraic equation in the corollary of Main Conjecture follows immediately 
from inserting $d$ as defined in Main Conjecture into equation (\ref{preliminaryCekva}) in Section 3 and letting $n\to\infty$.
 If we put $d$ into (\ref{preliminaryCekva}) we namely get
\begin{eqnarray}\label{Insertd}
\sum_{j=1}^{k}\bigg(\sum_{i=0}^{\deg Q_j}\alpha_{j,i}
\frac{z^{i}}{n^{\max_{j\in[j_0+1,k]}\big(\frac{j-j_0}{j-\deg Q_j}\big)(j-i)+j_0-j}}\bigg)C^j(z)=1.
\end{eqnarray}
Denote by $N_{j,i}$ the exponent of $n$ in (\ref{Insertd}) for given $j$ and 
$i$. Thus  
 \begin{eqnarray*}
N_{j,i}=\max_{j\in[j_0+1,k]}\bigg(\frac{j-j_0}{j-\deg Q_j}\bigg)(j-i)+j_0-j.
\end{eqnarray*}
The terms in (\ref{Insertd}) for which this exponent is positive tend to zero as $n\to\infty$.

First we consider $j$ for which $\deg Q_j=j$, and denote, as usual, by $j_0$ the largest such $j$. If $j=j_0$, then $i\leq \deg Q_{j_0}=j_0$ and thus for $j=j_0$ and $i=j_0$ we get
\begin{eqnarray*}
N_{j_0,j_0}&=&\max_{j\in[j_0+1,k]}\bigg(\frac{j-j_0}{j-\deg Q_j}\bigg)(j-i)+j_0-j\\
&=&\max_{j\in[j_0+1,k]}\bigg(\frac{j-j_0}{j-\deg Q_j}\bigg)(j_0-j_0)+j_0-j_0=0,
\end{eqnarray*}
and for $j=j_0$ and $i<j_0$ we have 
\begin{eqnarray*}
N_{j_0,i}&=&\max_{j\in[j_0+1,k]}\bigg(\frac{j-j_0}{j-\deg Q_j}\bigg)(j-i)+j_0-j\\
&>&\max_{j\in[j_0+1,k]}\bigg(\frac{j-j_0}{j-\deg Q_j}\bigg)(j_0-j_0)+j_0-j_0=0.
\end{eqnarray*}
Thus $N_{j_0,j_0}=0$ and $N_{j_0,i}>0$ for $i<j_0$, and for the term corresponding to $j=j_0$ in (\ref{Insertd}) we get 
\begin{eqnarray*}
\sum_{i=0}^{j_0}\alpha_{j_0,i}
\frac{z^{i}}{n^{\max_{j\in[j_0+1,k]}\big(\frac{j-j_0}{j-\deg Q_j}\big)(j_0-i)+j_0-j_0}}
C^{j_0}(z)
\to\alpha_{j_0,j_0}z^{j_0}C^{j_0}(z)=z^{j_0}C^{j_0}(z)
\end{eqnarray*}
 when $n\to\infty$ and using $\alpha_{j_0,j_0}=1$.\\
Now let $j$ be such that $\deg Q_j=j$ and $j<j_0$. Then 
$i\leq\deg Q_j=j$ and 
\begin{eqnarray*}
N_{j,j}&=&\max_{j\in[j_0+1,k]}\bigg(\frac{j-j_0}{j-\deg Q_j}\bigg)(j-i)+j_0-j\\
&=&\max_{j\in[j_0+1,k]}\bigg(\frac{j-j_0}{j-\deg Q_j}\bigg)(j-j)+j_0-j=j_0-j>0,
\end{eqnarray*}
and for $i<j$ we get 
\begin{eqnarray*}
N_{j,i}&=&\max_{j\in[j_0+1,k]}\bigg(\frac{j-j_0}{j-\deg Q_j}\bigg)(j-i)+j_0-j\\
&>&\max_{j\in[j_0+1,k]}\bigg(\frac{j-j_0}{j-\deg Q_j}\bigg)(j-j)+j_0-j=j_0-j>0,
\end{eqnarray*}
that is $N_{j,i}>0$ for all $j<j_0$ such that $\deg Q_j=j$ and for all $i\leq j$.
Thus for the corresponding terms in (\ref{Insertd}) we get
 \begin{eqnarray*}
\sum_{j\in\{j<j_0:\deg Q_j=j\}}
\sum_{i=0}^{\deg Q_{j}}\alpha_{j,i}
\frac{z^{i}}{n^{\max_{j\in[j_0+1,j]}\big(\frac{j-j_0}{j-\deg Q_j}\big)(j-i)+j_0-j}}
C^j(z)
\to 0
\end{eqnarray*}
when $n\to\infty$ for every $j<j_0$ for which $\deg Q_j=j$.\\ 

Now denote by $j_m$ the $j$ for which the maximum 
$d=\max_{j\in[j_0+1,k]}\big(\frac{j-j_0}{j-\deg Q_j}\big)$ is attained. Note that there may be several distinct $j$ for which this maximum is attained.\footnote{Consider for example the Laplace type operator (that is with all polynomial coefficients $Q_j$ linear) 
$T=zD+zD^2+\ldots zD^k$. Here $j_0=1$ and the equation satisfied by the Cauchy transform of the asymptotic root measure of the scaled eigenpolynomial $q_n(z)=p_n(nz)$ is given by $zC(z)+zC^2(z)+\ldots zC^k(z)=1$, since $d=\max_{j\in[2,k]}\big(\frac{j-j_0}{j-\deg Q_j}\big)=1$ is attained for \textit{every} $j=2,3,\ldots k$.}
Then 
\begin{eqnarray*}
N_{j_m,deg Q_{j_m}}&=&\max_{j\in[j_0+1,k]}\bigg(\frac{j-j_0}{j-\deg Q_j}\bigg)(j-i)+j_0-j\\
&=&\bigg(\frac{j_m-j_0}{j_m-\deg Q_{j_m}}\bigg)(j_m-\deg Q_{j_m})+j_0-j_m\\
&=&j_m-j_0+j_0-j_m=0,
\end{eqnarray*}
and for $i<\deg Q_{j_m}$ we get 
\begin{eqnarray*}
N_{j_m,i}&=&\max_{j\in[j_0+1,k]}\bigg(\frac{j-j_0}{j-\deg Q_j}\bigg)(j-i)+j_0-j\\
&>&\bigg(\frac{j_m-j_0}{j_m-\deg Q_{j_m}}\bigg)(j_m-\deg Q_{j_m})+j_0-j_m\\
&=&j_m-j_0+j_0-j_m=0,
\end{eqnarray*}
i.e. $N_{j_m,\deg Q_{j_m}}=0$ and $N_{j_m,i}>0$ for $i<\deg Q_{j_m}$, 
and for the term corresponding to $j=j_m$ in (\ref{Insertd}) we get 
\begin{eqnarray*}
\sum_{i=0}^{\deg Q_{j_m}}\alpha_{j_m,i}
\frac{z^{i}}{n^{\max_{j\in[j_0+1,k]}\big(\frac{j-j_0}{j-\deg Q_j}\big)(j_m-i)+j_0-j_m}}
C^{j_m}(z)\to \alpha_{j_m,\deg Q_{j_m}}z^{\deg Q_{j_m}}C^{j_m}(z)
\end{eqnarray*}
when $n\to\infty$. In case of several $j$ for which $d$ is attained, we put $A=\{j:(j-j_0)/(j-\deg Q_j)=d\}$, and for the corresponding terms in (\ref{Insertd}) we get 
\begin{eqnarray*}
\sum_{j\in A}
\sum_{i=0}^{\deg Q_{j}}\alpha_{j,i}
\frac{z^{i}}{n^{\max_{j\in[j_0+1,k]}\big(\frac{j-j_0}{j-\deg Q_j}\big)(j-i)+j_0-j}}
C^{j}(z)\to \sum_{j\in A}\alpha_{j,\deg Q_{j}}z^{\deg Q_j}C^{j}(z)
\end{eqnarray*}
when $n\to\infty$. 
Now consider the remaining terms in (\ref{Insertd}), namely terms for which $j<j_0$ such that
 $\deg Q_j<j$, terms for which $j_0<j<j_m$, and terms for which $j_m<j\leq k$ 
(clearly this last case does not exist if $j_m=k$).

 We start with $j<j_0$ such that $\deg Q_j<j$. Then  
$i\leq \deg Q_j<j$ and 
\begin{eqnarray*}
N_{j,i}&=&\max_{j\in[j_0+1,k]}\bigg(\frac{j-j_0}{j-\deg Q_j}\bigg)(j-i)+j_0-j\\
&>&
\max_{j\in[j_0+1,k]}\bigg(\frac{j-j_0}{j-\deg Q_j}\bigg)(j-j)+j_0-j=j_0-j>0,
\end{eqnarray*}
and for the corresponding terms in (\ref{Insertd}) we have 
\begin{eqnarray*}
\sum_{j\in\{j<j_0:\deg Q_j<j\}}\sum_{i=0}^{\deg Q_j}\alpha_{j,i}
\frac{z^{i}}{n^{\max_{j\in[j_0+1,k]}\big(\frac{j-j_0}{j-\deg Q_j}\big)(j-i)+j_0-j}}
C^j(z)\to 0
\end{eqnarray*}
when $n\to\infty$.\\

Now assume that $j_m<k$ and consider $j_m<j\leq k$. Clearly $j_m>j_0$ since the maximum is taken over $j\in[j_0+1,k]$, and therefore $i\leq \deg Q_j<j$ for
$j_m<j\leq k$. Also,  
\begin{eqnarray*}
\max_{j\in[j_0+1,k]}\bigg(\frac{j-j_0}{j-\deg Q_j}\bigg)=
\bigg(\frac{j_m-j_0}{j_m-\deg Q_{j_m}}\bigg)>
\bigg(\frac{j-j_0}{j-\deg Q_j}\bigg),
\end{eqnarray*}
since the maximum is attained for $j_m$ by assumption. Thus we get
\begin{eqnarray*}
N_{j,i} &=&\max_{j\in[j_0+1,k]}\bigg(\frac{j-j_0}{j-\deg Q_j}\bigg)(j-i)+j_0-j=
\bigg(\frac{j_m-j_0}{j_m-\deg Q_{j_m}}\bigg)(j-i)+j_0-j\\
&>& \bigg(\frac{j-j_0}{j-\deg Q_j}\bigg)(j-i)+j_0-j\geq  
\bigg(\frac{j-j_0}{j-\deg Q_j}\bigg)(j-\deg Q_j)+j_0-j\\
&=&j-j_0+j_0-j=0,
\end{eqnarray*}
i.e. $N_{j,i}>0$ for every $j_m<j\leq k$ and every $i\leq \deg Q_j$. For the corresponding terms in (\ref{Insertd}) we therefore get
\begin{eqnarray*}
\sum_{j_m<j\leq k}\sum_{i=0}^{\deg Q_j}\alpha_{j,i}
\frac{z^{i}}{n^{\max_{j\in[j_0+1,k]}\big(\frac{j-j_0}{j-\deg Q_j}\big)(j-i)+j_0-j}}
C^j(z)\to 0
\end{eqnarray*}
as $n\to\infty$.\\

 Finally we consider $j_0<j<j_m$. Note that this also covers the case $j_{m_1}<
j<j_{m_2}$ where the maximum $d$
is attained for both $j_{m_1}$ and $j_{m_2}$.
Since $i\leq \deg Q_j<j$ we get
\begin{eqnarray*}
N_{j,i}&=&\max_{j\in[j_0+1,k]}\bigg(\frac{j-j_0}{j-\deg Q_j}\bigg)(j-i)+j_0-j=
\bigg(\frac{j_m-j_0}{j_m-\deg Q_{j_m}}\bigg)(j-i)+j_0-j\\
&>&\bigg(\frac{j-j_0}{j-\deg Q_{j}}\bigg)(j-i)+j_0-j\geq 
\bigg(\frac{j-j_0}{j-\deg Q_{j}}\bigg)(j-\deg Q_j)+j_0-j\\
&=&
j-j_0+j_0-j=0,
\end{eqnarray*}
i.e. $N_{j,i}>0$ for every $j_0<j<j_m$ and every $i\leq \deg Q_j$. Thus for the corresponding terms in (\ref{Insertd}) we get
\begin{eqnarray*}
\sum_{j_0<j<j_m}\sum_{i=0}^{\deg Q_j}\alpha_{j,i}
\frac{z^{i}}{n^{\max_{j\in[j_=+1,k]0}\big(\frac{j-j_0}{j-\deg Q_j}\big)(j-i)+j_0-j}}C^j(z)\to 0
\end{eqnarray*}
when $n\to\infty$.\\\\
Adding up these results we finally get the following equation by letting $n\to\infty$
in equation (\ref{Insertd}):  
\begin{displaymath}
z^{j_0}C^{j_0}(z)+ \sum_{j\in A}\alpha_{j,\deg Q_{j}}z^{\deg Q_j}
C^j(z)=1,
\end{displaymath}
where $j_0$ is the largest $j$ for which $\deg Q_j=j$, and $A$ is the set consisting of all $j$ for which the maximum $d=\max_{j\in [j_0+1,k]}\big(\frac{j-j_0}{j-\deg Q_j}\big)$ is attained, i.e. 
$A=\{j:(j-j_0)/(j-\deg Q_j)=d\}$.
\subsection{Theorem 5}
Here we prove that for a class of operators containing the operators considered in Corollaries 1 and 2 the conjectured upper bound 
$\lim_{n\to\infty}\sup(r_n/n^d)\leq c_1$ implies the conjectured 
lower bound 
$\lim_{n\to\infty}\inf(r_n/n^d)\geq c_0$ for some constants $c_1\geq c_0>0$ and where $d$ is as in Main Conjecture. This follows automatically from inequality (\ref{lemma3}) in Lemma 3. We have the following\\\\
\textbf{Theorem 5.} \textit{Let $T$ be a degenerate exactly-solvable operator of order $k$ which satisfies the condition}
\begin{eqnarray*}
b:=\min_{j\in[1,k-1]}^{+}\bigg(\frac{k-j}{k-j+\deg Q_j-\deg Q_k}\bigg)
=\max_{j\in[j_0+1,k]}\bigg(\frac{j-j_0}{j-\deg Q_j}\bigg)=:d,
\end {eqnarray*}
\textit{where the notation 
 $\min^{+}$ means that the minimum is taken only over 
positive values of $(k-j+\deg Q_j-\deg Q_k)$.
Assume that the inequality $r_n\leq c_1(n-k+1)^{d}$ holds for some positive constant $c_1$ for all sufficiently large $n$.
Then there exists a positive constant $c_0\leq c_1$ such that 
$r_n\geq c_0(n-k+1)^{d}$ holds for all sufficiently large $n$. Thus}
\begin{displaymath}
\lim_{n\to\infty}\sup\frac{r_n}{n^d}\leq c_1\quad\Rightarrow\quad
 \lim_{n\to\infty}\inf\frac{r_n}{n^d}\geq c_0.
\end{displaymath}
\\
\textbf{Proof.} 
From inequality (\ref{lemma3}) in Lemma 3 we have 
\begin{eqnarray}\label{appendix}
1&\leq& \sum_{j=1}^{k-1}\sum_{i=0}^{\deg Q_j}|\alpha_{j,i}|2^{k-j}\frac{r_n^{k-j+i-\deg Q_k}}{(n-k+1)^{k-j}}+\sum_{0\leq i<\deg Q_k}\frac{|\alpha_{k,i}|}{r_{n}^{\deg Q_k-i}}\nonumber\\
&\leq & \sum_{j=1}^{k-1}K_j\frac{r_n^{k-j+\deg Q_j-\deg Q_k}}{(n-k+1)^{k-j}}
+\sum_{0\leq i<i_k}\frac{|\alpha_{k,i}|}{r_{n}^{i_k-i}}
\end{eqnarray}
where the $K_j$ are positive constants. The second sum on the right-hand side of (\ref{appendix}) tends to zero when $n\to\infty$ due to Theorem 1. 
We now decompose the first sum on the right-hand side of (\ref{appendix}) into three parts. Namely, let\\\\
$A=\{j:\frac{k-j}{k-j+\deg Q_j-\deg Q_k}=d\}$, and note that $(k-j+\deg Q_j-\deg Q_k)>0$ here since $d>0$.\\\\
$B=\{j:\frac{k-j}{k-j+\deg Q_j-\deg Q_k}> d\}$, and note that $(k-j+\deg Q_j-\deg Q_k)>0$ here since $d>0$.\\\\
$C=\{j:(k-j+\deg Q_j-\deg Q_k)\leq 0\}$, and note that $j<k$ in  (\ref{appendix}).\\\\
Clearly due to the condition $b=d$ there are no terms for which 
 $\frac{k-j}{k-j+\deg Q_j-\deg Q_k}<d$ and 
$(k-j+\deg Q_j-\deg Q_k)>0$ both hold .\\\\ 
If $j\in A$ then 
\begin{displaymath} 
\frac{r_n^{k-j+\deg Q_k-\deg Q_k}}{(n-k+1)^{k-j}}=
\bigg(\frac{r_n}{(n-k+1)^{d}}\bigg)^{k-j+\deg Q_j-\deg Q_k}
\end{displaymath}
for the corresponding terms in the sum on the right-hand side of 
(\ref{appendix}).\\\\
If $j\in B$ then 
$d(k-j+\deg Q_j-\deg Q_k)<(k-j)$, and
  this inequality together with the upper bound 
$r_n\leq c_1(n-k+1)^{d}$ which we assume holds for all sufficiently large $n$, gives us 
\begin{displaymath}
\frac{r_n^{k-j+\deg Q_j-\deg Q_k}}{(n-k+1)^{k-j}}\leq 
\frac{c_1(n-k+1)^{d(k-j+\deg Q_j-\deg Q_k)}}
{(n-k+1)^{k-j}}\to 0 
\end{displaymath}
when $n\to\infty$ for
the corresponding terms in (\ref{appendix}).\\\\
If $j\in C$ then $(k-j+\deg Q_j-\deg Q_k)\leq 0$ and we get
\begin{displaymath}
\frac{r_n^{k-j+\deg Q_j-\deg Q_k}}{(n-k+1)^{k-j}}\to 0
\end{displaymath} 
when $n\to\infty$ 
for the corresponding terms in (\ref{appendix}) due to Theorem 1. Note that if $(k-j+\deg Q_j-\deg Q_k)=0$ the corresponding term tends to zero when $n\to\infty$ since $j<k$ in (\ref{appendix}).\\\\
With this decomposition of the first sum on the right-hand side of the last inequality in (\ref{appendix}) we can write
\begin{eqnarray*}
1&\leq& 
\sum_{j=1}^{k-1}K_j\frac{r_n^{k-j+\deg Q_j-\deg Q_k}}{(n-k+1)^{k-j}}
+\sum_{0\leq i<i_k}\frac{|\alpha_{k,i}|}{r_{n}^{i_k-i}}\\
&\leq & \sum_{j\in A}K_j\bigg(\frac{r_n}{(n-k+1)^{d}}\bigg)^{k-j+\deg Q_j-\deg Q_k}\\
&+&\sum_{j\in B}
K_j\frac{c_1(n-k+1)^{d(k-j+\deg Q_j-\deg Q_k)}}{(n-k+1)^{k-j}}\\
&+&\sum_{j\in C}^{}K_j\frac{r_n^{k-j+\deg Q_j-\deg Q_k}}{(n-k+1)^{k-j}} 
+\sum_{0\leq i<i_k}\frac{|\alpha_{k,i}|}{r_{n}^{i_k-i}}
\end{eqnarray*}
where the last three sums tend to zero when $n\to\infty$ by the above arguments (and the last one due to Theorem 1).\\\\ 
Thus for all sufficiently large $n$ there exists a positive 
constant $c'$ such that 
\begin{eqnarray}\label{hej}
c'\leq \sum_{j\in A}K_j\bigg(\frac{r_n}{(n-k+1)^{d}}\bigg)^{k-j+\deg Q_j-\deg Q_k}
\end{eqnarray}
where $c'\to 1$ when $n\to\infty$. If the set $A$ contains precisely one element, then the sum in (\ref{hej}) consists of one single term, and we are done: there exists a positive constant $c_0=(c'/K_j)^{1/(k-j+\deg Q_j-\deg Q_k)}$ such that $r_n\geq c_0(n-k+1)^{d}$ for all sufficiently large $n$, and thus 
$\lim_{n\to\infty}\inf(r_n/n^d)\geq c_0$. 

But
clearly, for some operators $A$ will contain more than one element.
If this is the case define $m:=\min_{j\in A}(k-j+\deg Q_j-\deg Q_k)$
and denote by $j_m$ the corresponding $j$ for which this minimum
 is attained. Using the upper bound $r_n\leq c_1(n-k+1)^{d}$ we then get the following inequality\footnote{Consider for example the operator $T=zD+D^2+zD^3+zD^4$. Here by Lemma 3:
\begin{eqnarray*}
1\leq\sum_{j=1}^{3}\frac{2^{4-j}r_n^{3-j+\deg Q_j}}{(n-3)^{4-j}}= 
8\frac{r_n^3}{(n-3)^3}+4\frac{r_n}{(n-3)^2}+2\frac{r_n}{(n-3)}, 
\end{eqnarray*}
where $r_n$ is the largest modulus of all roots of the unique and monic eigenpolynomial of $T$. For this operator $d=1$ and we see that 
$(4-j)/(3-j+\deg Q_j)=d$ for $j=1$ and for $j=3$.
Now, assuming that $r_n\leq c_1(n-3)$ holds for some positive constant $c_1$ for large $n$, our inequality becomes
\begin{eqnarray*}
1&\leq& 8\frac{r_n^3}{(n-3)^3}+4\frac{r_n}{(n-3)^2}+2\frac{r_n}{(n-3)}\\
&\leq&
8\frac{r_n}{(n-3)}\cdot\frac{c_1^2(n-3)^2}{(n-3)^2}+4\frac{c_1(n-3)}{(n-3)^2}
+2\frac{r_n}{(n-3)}\\
&=&(8c_1^2+2)\frac{r_n}{(n-3)}+\frac{4c_1}{(n-3)}
\end{eqnarray*}
where the last term tends to zero as $n\to\infty$. Thus
$r_n\geq c_0(n-3)$ 
for sufficiently large choices on $n$, where $c_0=1/(8c^2_1+2)$, and hence 
$\lim_{n\to\infty}\inf(r_n/n)\geq c_0$.} from (\ref{hej}):
\begin{eqnarray*}
c'&\leq& \sum_{j\in A}K_j\bigg(\frac{r_n}{(n-k+1)^{d}}\bigg)^{k-j+\deg Q_j-\deg Q_k}\\
&=& K_{j_m}\bigg(\frac{r_n}{(n-k+1)^{d}}\bigg)^m\\
&+&
\sum_{j\in A\backslash\{j_m\}}K_j\bigg(\frac{r_n}{(n-k+1)^{d}}\bigg)^m\cdot\bigg(\frac{r_n}{(n-k+1)^{d}}\bigg)^{k-j+\deg Q_j-\deg Q_k-m}\leq 
\end{eqnarray*}
\begin{eqnarray*}
&\leq &  K_{j_m}\bigg(\frac{r_n}{(n-k+1)^{d}}\bigg)^m+
\sum_{j\in A\backslash \{j_m\}}K_j\bigg(\frac{r_n}{(n-k+1)^{d}}\bigg)^m\cdot c_1^{k-j+\deg Q_j-\deg Q_k-m}\\
&=&
\bigg(\frac{r_n}{(n-k+1)^{d}}\bigg)^m\bigg(K_{j_m}+\sum_{j\in A\backslash \{j_m\}}K_j\cdot c_1^{k-j+\deg Q_j-\deg Q_k-m}\bigg)\\
&=&
\bigg(\frac{r_n}{(n-k+1)^{d}}\bigg)^m\cdot K
\end{eqnarray*}
where $K>0$. Thus $r_n\geq \big(\frac{c'}{K}\big)^{1/m}(n-k+1)^d$ for all sufficiently large $n$, and therefore there exists a positive constant $c_0=(1/K)^{1/m}$ (recall that $c'\to 1$ when $n\to\infty$) such that $\lim_{n\to\infty}\inf\frac{r_n}{n^d}\geq c_0$.
\hfill$\square$\\\\\\\\
\section{Open Problems}
\textbf{1.} The main challenge is to obtain a complete proof of the Main Conjecture, see Introduction. This proof requires both the sharp upper and lower bounds of the largest root. The upper bound can apparently be obtained by a detailed study of the corresponding Riccati equation at $\infty$. If the Main Conjecture is settled then to achieve its corollary (the Cauchy transform equation) one can use a technique similar to that of \cite{BR} to prove the basic assumption, see Section 3.\\\\
\textbf{2.} As suggested by one of the referees, estimates similar to that of Main Conjecture can be formulated for the sequence of roots $z_{n,i}$ of $p_n$ such that $\lim_{n\to\infty}\frac{r_{n,i}}{r_{n,n}}=\alpha$ where $0<\alpha<1$ and  $|z_{n,i}|=r_{n,i}$.\\\\
\textbf{3.} In the case of orthogonal polynomials there is a number of results describing the growth of the largest modulus $r_n$ of the roots as an expansion in powers of $n$, see e.g. \cite{KI1} and \cite{KI2} and references therein. In the present paper we conjectured the form of the leading term of $r_n$ in our more general setup. As suggested by one of the referees, the question about the lower terms in the expansion of $r_n$ is natural in our context as well.\\\\
\textbf{4.} Operators of the type we consider occur in the theory of Bochner-Krall orthogonal systems, i.e. families of polynomials which are both eigenfunctions of some finite order differential operator and orthogonal with respect to some suitable inner product. A lot is known about the asymptotic zero 
distribution of orthogonal polynomials, and by comparing such known results with results on the asymptotic zero distribution of eigenpolynomials of degenerate exactly-solvable operators, we believe it will be possible to gain new insight into the nature of BKS. We have previously used our results from \cite{BR} to prove a special case of a general conjecture describing the leading terms of all Bochner-Krall operators, see \cite{BRS}. 
Another problem relevant for BKS is to describe all exactly-solvable operators whose eigenpolynomials have real roots only.\\\\
\textbf{5.} Numerical evidence indicates that the roots of the scaled eigenpolynomials fill certain curves in the complex plane. The support of the limiting 
root measure $\mu$ seems to be a tree. This is the case for the non-degenerate exactly-solvable operators which we treated in \cite{BR}, but then without such a scaling of the eigenpolynomials. By a tree we mean a connected compact subset $\Gamma$ of $\mathbb{C}$ which consists of a finite union of analytis curves and where $\mathbb{\hat{C}}\setminus \Gamma$ is simply connected. 
The (conjectural) algebraic equation satisfied by the Cauchy transform contains a lot of information about $\mu$, and it remains to describe its support explicitly.  \\\\
\textbf{6.} Conjecturally
the support of the asymptotic zero distribution of the scaled eigenpolynomial $q_n$ is the union of a finite number of analytic curves in the complex plane which we denote by $\Xi_T$, i.e. $\Xi_T= $ supp $\mu$, where $\mu$ is the limiting root measure of $q_n$. Then the following conjecture seems to be quite plausible\footnote{The question concerning interlacing was raised by B. Shapiro. Also see \cite{Bender}.}.\\\\
\textbf{Conjecture 1.} [Interlacing property] \textit{For any family $\{q_n\}$ of appropriately scaled eigenpolynomials of a degenerate exactly-solvable operator, the zeros of any two consecutive polynomials $q_{n+1}$ and $q_n$ interlace along $\Xi_T$ for all sufficiently large integers $n$.}\\

When defining the interlacing property some caution is required since the zeros of $q_n$ do not lie exactly on $\Xi_T$. Thus identify some sufficiently small neighbourhood $N(\Xi_T)$ of $\Xi_T$ with the normal bundle to $\Xi_T$ by equipping $N(\Xi_T)$ with the projection onto $\Xi_T$ along the fibres which are small curvilinear segments orthogonal to $\Xi_T$. We then say that two sets of points in $N(\Xi_T)$ interlace if their orthogonal projections on $\Xi_T$ interlace in the usual sense. If $\Xi_T$ has singularities one should first remove some sufficiently small neighbourhoods of these singularities and then proceed as above on the remaining part of $\Xi_T$. Conjecture 1 thus states that for any sufficiently small neighbourhood $N(\Xi_T)$ of $\Xi_T$ there exists a number $n_0$ such that the interlacing property holds for the zeros of $q_n$ and $q_{n+1}$ for all $n\geq n_0$. We conclude this section by showing some pictures illustrating the interlacing property. 
Below, small dots represent the roots of $q_{n+1}$ and large dots represent the roots of $q_n$ for some fixed $n$.\\

\begin{tabular}{cc}
\includegraphics[width=4.0cm]{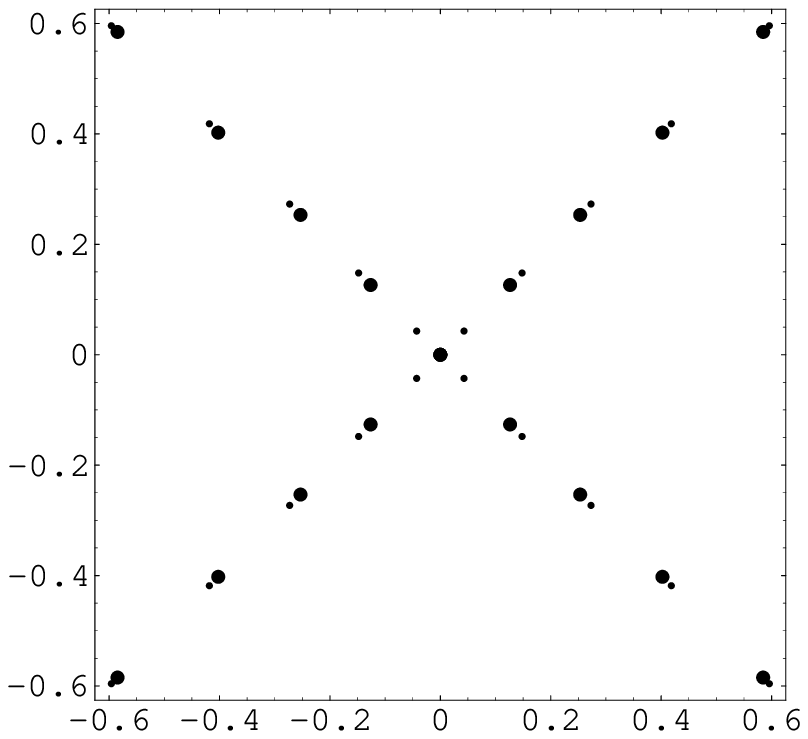} & 
\includegraphics[width=4.0cm]{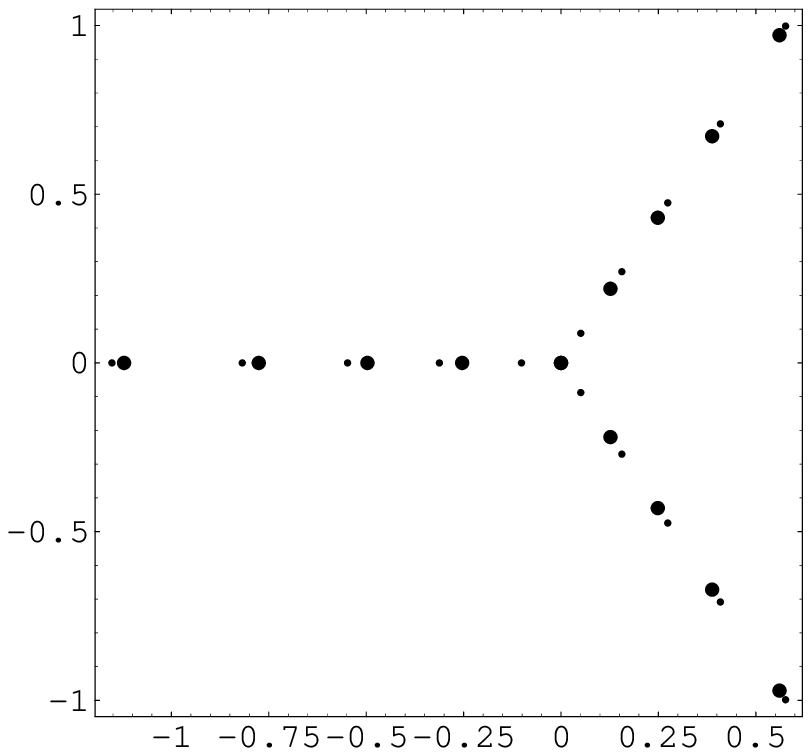}\\
$\qquad T=z^2D^2+z^3D^3+zD^5\qquad$, &  $\qquad T=zD+z^2D^2+D^3\qquad$,\\
roots of $q_{25}$ and $q_{24}$. & roots of $q_{20}$ and $q_{19}$.
\end{tabular}\\

\begin{tabular}{cc}
\includegraphics[width=4.0cm]{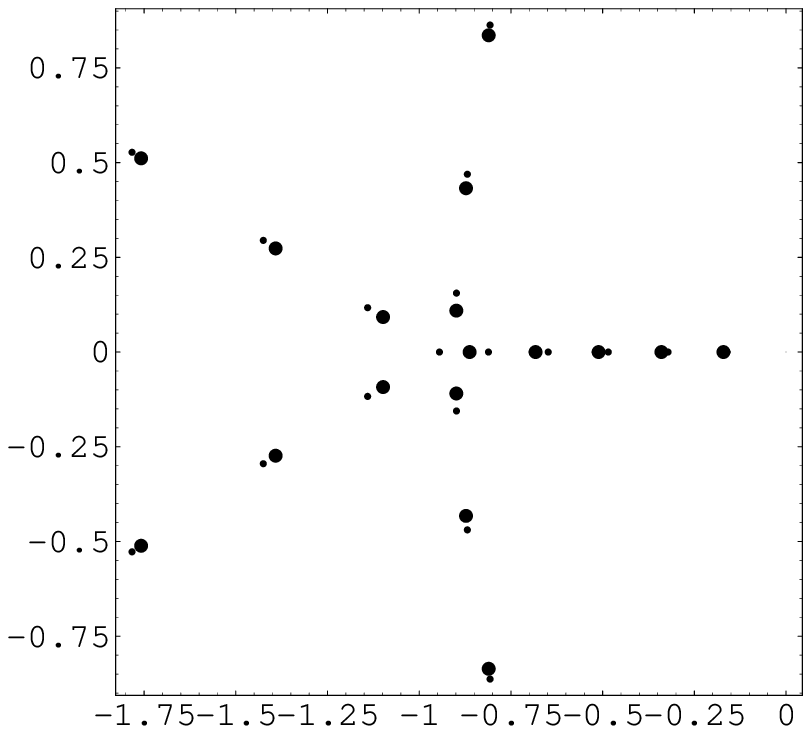} &
\includegraphics[width=4.0cm]{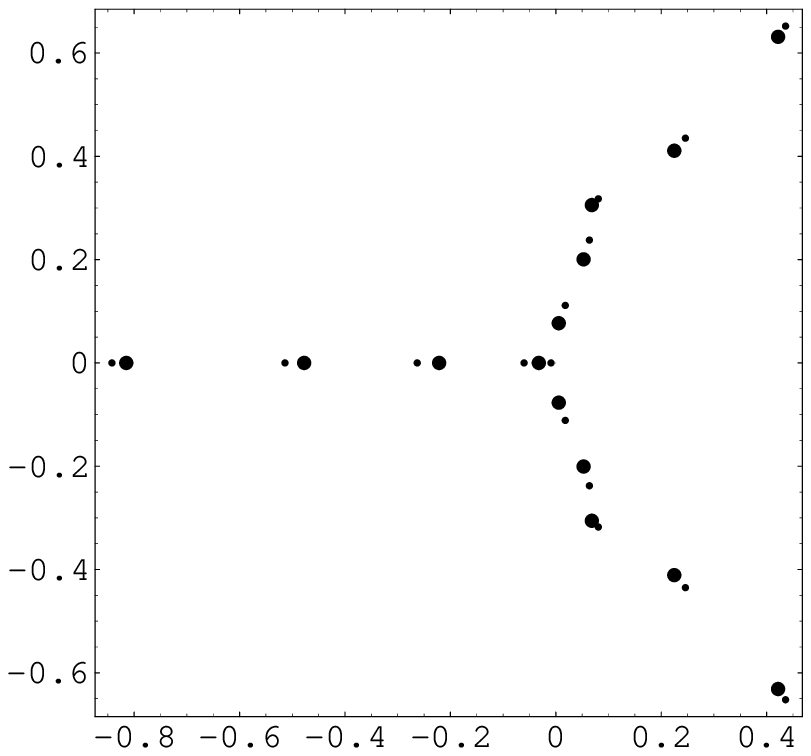}\\
$T=zD+zD^2+zD^3+zD^4+zD^5$, & $T=z^3D^3+z^2D^5+zD^6$,\\
roots of $q_{23}$ and $q_{22}$. & roots of $q_{20}$ and $q_{19}$.
\end{tabular}\\

\end{document}